\documentclass[11pt, a4paper]{article}
\usepackage[subpreambles=false]{standalone}%
\usepackage[english]{babel}
\usepackage[top=1.5cm, bottom=1.5cm, left=1.5cm, right=1.5cm, headheight=.5cm, headsep=.5cm, a4paper]{geometry}
\usepackage{caption}
\usepackage{subcaption}
\usepackage{float}
\usepackage{xcolor}
\usepackage{graphicx}
\usepackage{graphics}
\usepackage[unicode, colorlinks = true]{hyperref}
\hypersetup{
colorlinks=true,
linkcolor=[rgb]{.75,.15,.10},
citecolor=[rgb]{.35,.65,.55},
urlcolor=[rgb]{.88,.33,.88}}

\usepackage[shortlabels]{enumitem}
\usepackage[parfill]{parskip}
\usepackage{amsmath,amsthm,amssymb,amsfonts}
\usepackage{mathbbol}
\usepackage{dsfont}
\usepackage{pgfplots}
\usepackage[normalem]{ulem}
\pgfplotsset{compat=1.15}
\usepackage{mathrsfs}
\usetikzlibrary{arrows}
\usepackage{quiver}
\usepackage{authblk}

\def\d{{\mathrm d}}

\def\L{\mathcal L}

\newcommand\restr[2]{
\left.#1 \vphantom{\Bigr|}\right\lvert_{#2}}

\def\abs#1{\left\lvert#1\right\rvert}

\def\1{\mathds 1}
\def\0{\varnothing}

\def\wt{\widetilde}

\DeclareMathOperator{\sech}{sech}

\begin{document}
\sloppy
\def\binoppenalty{10000}
\def\phi{\varphi}
\sloppy

\title{Reconstruction of the non-linear wave at a buoy from shoreline data and applications to the tsunami inverse problem for piece-wise sloping bathymetry\footnote{To be published in Water Waves}}
\author[1]{Oleksandr Bobrovnikov\footnote{Corresponding author (obobrovnikov@alaska.edu)}}
\author[2]{Madison Jones}
\author[3]{Shriya Prasanna}
\author[2]{Josiah Smith}
\author[1]{Alexei Rybkin}
\author[4,5]{Efim Pelinovsky}

\affil[1]{University of Alaska Fairbanks, Fairbanks, Alaska, United States}
\affil[2]{University of Colorado Boulder, Boulder, Colorado, United States}
\affil[3]{University of Washington, Seattle, Washington, United States}
\affil[4]{HSE University, Nizhny Novgorod, Russia}
\affil[5]{Institute of Applied Physics, Nizhny Novgorod, Russia}
\maketitle

\begin{abstract}
We discuss the following inverse problem: given the run-up data of a tsunami wave, can we recover its initial shape? We study this problem within the framework of the non-linear shallow water equations, a model widely used to study tsunami propagation and inundation. Previously, we demonstrated that in the case of infinite sloping bathymetry, it is possible to recover the initial water displacement and velocity from shoreline readings \cite{Rybkin23,Rybkin24,Rybkin25}. 

We consider a finite sloping bathymetry. We show that it is possible to recover boundary conditions (water displacement and velocity) on a virtual buoy from the shoreline data. Further, we discuss stitching together the shallow water equations and the Boussinesq equation in a more complex piece-wise sloping bathymetry in order to recover the 
initial conditions, while incorporating the dispersion into our model. 
\end{abstract}


\section*{Acknowledgments}
Part of this work was done during the 2025 summer REU program run by Dr.~Alexei~Rybkin and was supported by National Science Foundation grant DMS-2307774.
OB, MJ, SP, JS, and AR acknowledge support from NSF grant DMS-2307774.
The work of EP is supported by the project of Mathematical Centre ``Symmetry, Information, Chaos" within the Program of fundamental research of HSE in 2025 (section 1).
We also thank the UAF DMS for hosting us.

\section{Introduction}\label{intro}
In 1964, the 9.2\,$M_w$ Good Friday earthquake and tsunami resulted in 124 fatalities and millions of dollars in economic devastation \cite{Wood2014}. This event, as well as the 1960 Valdivia, the 2004 Indian Ocean, and the 2011 Tōhoku tsunamis, demonstrates the significant threats to life and infrastructure that tsunamis pose to coastal communities.

Thus, there is a strong incentive to understand this phenomenon to improve mitigation strategies. Mathematical modelling plays a key role, allowing for the prediction of wave behaviour and the development of warning systems or preventative architecture such as seawalls. We suggest the works of \cite{Dias07, Dias14, levin2016} for an overview of the lifespan of a tsunami wave: generation, propagation, and inundation.

Tsunamis are often studied within 2+1 (2 spatial and 1 temporal variable) models, such as non-linear shallow water equations (NSWE), the Boussinesq equation, the Korteweg-De Vries equation, etc. They are obtained from Euler equations, a non-linear 3+1 system, utilising various approximations such as depth-averaging, long wave assumption, zero transverse variations, assumption of no vorticity, dispersion, etc. \cite{Johnson97}.
It is debatable whether the dispersion matters when modelling tsunami propagation. Some argue that under certain circumstances, depending on travel distance, mode of generation, and size of the wave, dispersive effects are noticeable \cite{DispersionPro}, while others argue that they do not matter \cite{DispersionCon}.

There appear to be four ways to approach the inversion of tsunamis. The first method involves combining seismic and waveform data \cite{SeismicInv1, SeismicInv2}. This approach is susceptible to underestimating the incoming wave, but it relies on easily measurable data, such as teleseismic, geodetic, tide gauges, and sea surface heights, which is available before the wave event occurs. The second approach is reconstructing the tsunami from only waveform data gathered from mareographs or combined with other non-seismic information \cite{WaveFormInv, TsunamiWaveformAdv}. The next approach utilises sediment displacement and distribution to approximate tsunami characteristics and invert the wave \cite{DepositInv1, DepositInv2}. The last approach uses perhaps the most important characteristic of tsunamis: their inundation \cite{ComparisonPaper, RunupInv}. For a more in-depth discussion on numerical inversion techniques, please refer to \cite{Satake21}

Our paper takes the latter approach. In the first part of this text, we present an analytical algorithm to recover the wave field on a (virtual) buoy given the run-up records. With the advent of the DART (Deep-ocean Assessment and Reporting of Tsunamis) system, many models can predict run-up and destruction with the data they provide \cite{DART, percival18, mori_etal}. Our project aims to supply a theoretical framework to inform buoy placement. Given records of either historical or synthetic run-up data, our model can provide the travel time and the shape of the approaching wave. This will help decide where to place the buoy. Close to the shore, the non-linear effects of the tsunami wave become noticeable. We therefore use the non-linear shallow water equations (NSWE), which are due to Saint-Venant. We use the Carrier-Greenspan transform, which offers exact linearisation of NSWE for a certain class of bathymetries; however it can only be used for non-breaking waves.

The second part of our paper suggests two extensions of our main result. We consider a piece-wise sloping region and stitch together NSWE on the finite sloping region and another model on the semi-infinite flat bottom. For the other model, we take the Boussinesq equation and the linear shallow water equation. 
A similar piece-wise sloping plane beach configuration was considered by \cite{Synolakis87}, where linear and non-linear shallow water equations were used to model wave propagation in both sloping and flat bottom regions. Other works that discuss stitching NSWE with Boussinesq-type equations are those of \cite{VARSOLIWALA2021367, rigal25}. Unlike our work, all of these papers seek to solve the direct (forward in time) problem.  

The rest of the paper is organised as follows: Section \ref{sec:finite_int} contains our main result, the solution to the inverse problem for the NSWE system on a finitely sloping interval. We use the Carrier-Greenspan transform to recover boundary conditions in the hodograph plane and then present a way to recover the wave field at a buoy away from the shoreline (boundary conditions in the physical space). In Section \ref{subsec:num}, we examine the validity of our model with physically realistic scenarios and recover expected results. In Section \ref{sec:bouss}, we discuss stitching the NSWE system in a sloping region with a Boussinesq model in a semi-infinite, flat-bottomed region. We explore a heuristic understanding of the model in \ref{subsec:bouss_heur} and perform the implementation using the assumption of a two-soliton solution in Subsection \ref{subsec:bouss_implementation}. In Section \ref{sec:compare}, we present yet another way to treat piece-wise sloping beach following work of \cite{Synolakis87}, and compare the two inversion models in the piece-wise sloping beach. Finally, we provide concluding remarks in Section \ref{sec:conclusions}, as well as some possible extensions of our work.

\section{Inverse Problem for NSWE on the Finite Interval} \label{sec:finite_int}
\begin{figure}[b]
   \begin{center}
    \begin{subfigure}{.45\textwidth}
     \input{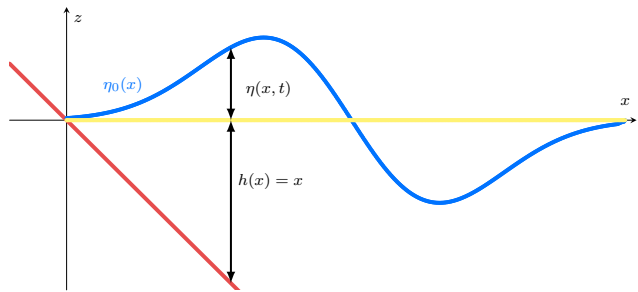}
    \end{subfigure}
   \end{center}
   \caption{Sketch of the bay geometry. Cross-sectional view; bathymetry is in red, the unperturbed water level is in yellow, and the water level is in blue.
   The total perturbed water depth is given by $H(x,t)=h(x)+\eta(x,t)$.}\label{fig:bath}
\end{figure}

We consider a finite, power-shaped, sloping beach of length $L$ as depicted in Figure \ref{fig:bath}, governed by the non-linear shallow water equations system (NSWE)
\begin{equation}\label{eq:SWE}
   \begin{aligned}
      &\partial_t\eta+u \partial_x (x+\eta) + (x+\eta)\partial_x u=0,&\text{(mass)}  \\
      &\partial_t u+u\partial_x u+\partial_x \eta= 0. &\text{(momentum)}
    \end{aligned}
\end{equation}
Here \eqref{eq:SWE} is given in dimensionless units: \(\eta(x, t)\) is the water displacement and \(u(x, t)\) is the velocity averaged over the vertical cross-section.
The substitution
\begin{equation}\label{eq:dimensions}
    \wt x = (H_0 / \alpha)  x, \quad
    \wt t = \sqrt{H_0 / g }\, t / \alpha,\quad
    \wt \eta = H_0  \eta,\quad
    \wt u = \sqrt{H_0 g} \, u,
\end{equation}
where $H_0 $ is a typical (characteristic) height, \(\alpha\) is the slope of the bay, and \(g\) is acceleration due to gravity, brings system \eqref{eq:SWE} to dimensional (with tildes) units.

This model is widely used to study tsunami propagation and inundation. There are multiple approaches to this problem: analytical, like those of \cite{Carrier58,Kanoglu06,Shimozono,Rybkin14}; numerical that use finite differences, elements, or volumes, like the works of \cite{KOUNADIS2020112315,XING2017361,Volna,MovingGrids}; or analytical-numerical \cite{NumAnExp,Bueler22}. More recent works use physics-informed neural networks for tsunami inundation modelling, as done in \cite{BRECHT2025114066}.
The works cited above are concerned with what we call the forward problem of tsunami modelling: from known initial or boundary conditions, one wants to find the run-up of the wave, or the motion of the free dry/wet boundary. In this subsection, we present a complete analytical solution to the following inverse problem: suppose that at zero time there is no water disturbance on the interval \([0,L]\); the goal is to recover the boundary conditions at \(x = L\) from observations at the shore. This part of our work is an extension of our previous results \cite{Rybkin23,Rybkin24,Rybkin25}, where we recovered the initial conditions from the shoreline oscillations in the case of infinite bathymetry in the \(x\) direction. System \eqref{eq:SWE} does not contain a dispersive term. In practice, dispersion can only be ignored when the wave is close to the shoreline, so our contribution here allows us to solve a more realistic problem.

We consider system \eqref{eq:SWE} for \((x, t) \in [0,L]\times [0,+\infty)\)  together with zero initial conditions
\begin{equation}\label{icond}
\begin{aligned}
     &u(x,0) = 0,  \\
     &\eta(x,0) = 0,
\end{aligned}
\end{equation}
and boundary conditions which imitate the wave field on the buoy
\begin{equation}
\begin{aligned}
     &u(L,t) = u_{\text{b}}(t),  \\
     &\eta(L,t) = \eta_{\text{b}}(t).
\end{aligned}
\end{equation}
Our goal is to recover boundary conditions from the vertical oscillations of the shoreline, run-up data $R (t)$.
We note here that the boundary conditions \( u_{\text{b}}(t), \eta_{\text{b}}(t)\) are assumed not to be assigned independently. See works of \cite{Antuono07,Petcu_Temam} for the discussion of the compatibility conditions for the boundary data. Remarkably, once compatibility conditions are imposed, we are able to recover two functions \(u(L,t)\) and \(\eta(L, t)\) from only one function \(R(t)\).

We utilise the Carrier-Greenspan transform (CGT) introduced by \cite{Carrier58} for the sloping plane beach. This transformation linearises \eqref{eq:SWE} and fixes the free boundary. It reads as
\begin{equation}\label{CG}
    \varphi(\sigma,\tau) = u(x,t),\quad \sigma = x + \eta(x,t),\quad \psi(\sigma,\tau) = \eta(x,t) + u^2(x,t)/2,\quad \tau = t-u(x,t).
\end{equation}  
Here \(\sigma, \tau\) are new variables and \(\psi, \phi\) are new functions.
System \eqref{eq:SWE} then under \eqref{CG}
transforms to
\begin{equation}\label{linearsys-pre}
    \begin{aligned}
    &\partial_\tau \psi + \sigma \partial_\sigma\varphi + \varphi = 0,\\
    &\partial_\tau\varphi + \partial_\sigma \psi = 0,
    \end{aligned}
\end{equation}
or
\begin{equation} \label{linearsys}
    \partial^2_\tau\psi = \sigma\partial^2_\sigma\psi + \partial_\sigma\psi.
\end{equation}
The cost of linearising \eqref{eq:SWE} is a complication of the initial and boundary conditions in the hodograph \((\sigma, \tau)\) plane. Additionally, the CGT only works as long as the wave does not break.
Since we consider zero initial conditions, there are no complications yet.
Initial conditions \eqref{icond} transform to
\begin{equation}
\begin{aligned}
     & \psi(\sigma,0) = 0,\\ 
     & \phi(\sigma, 0) = 0.
\end{aligned}
\end{equation}
In contrast, the boundary conditions are now specified on a parametric curve in the \((\sigma, \tau)\) plane, and the curve itself depends on the boundary conditions. Indeed, if we plug \(x = L\) into \eqref{CG}, we observe that in the second equation \(\sigma\) is not constant unless \(\eta(L,t) \equiv 0\). The method of data projections introduced by \cite{Rybkin21} allows one to replace boundary conditions on the parametric curve with equivalent conditions on the straight line \(\sigma = \sigma_L\), so that the solutions to the original and modified initial boundary value problem (IBVP) are within a desirable error tolerance. As noted previously in \cite{Rybkin25}, such complications are absent for the inverse problem, so we do not expand on this method. 

Once the compatibility condition and the method of data projections are brought into play, one obtains the IBVP 
\begin{equation}\label{eq:IBVP_psi}
    \begin{aligned}
       &\partial^2_\tau \psi = \sigma \partial^2_ \sigma \psi + \partial_\sigma \psi,\\ 
       & \psi(\sigma, 0) = 0, \\ 
       & \psi_\tau (\sigma, 0) = 0, \\ 
       & \psi(\sigma_L, \tau) = \psi_{\text b}(\tau),\\
       & \abs{\psi(0, \tau)} < \infty
    \end{aligned}
\end{equation}
on \((\sigma, \tau) \in [0,\sigma_L]\times [0,+\infty)\), which we refer to as the forward problem. For what follows, we assume that \(\psi_{\text{b}}\) satisfies \(\psi_{\text{b}}(0) = \psi_{\text{b}}'(0) = 0\). We briefly walk the reader through the solution to the direct problem provided by \cite{Rybkin21}.
We introduce the new variable $\rho^2 = \sigma /\sigma_L$ and the new function \({\theta(\rho, \tau) =  \psi(\sigma, \tau) - \psi_{\text{b}}(\tau)}\). This substitution homogenises the Dirichlet's boundary condition on the right end and reduces \eqref{linearsys} to the form
\begin{equation}\label{bessel}
   \partial^2_\tau \theta = \frac{1}{4 \sigma_L} \left( \partial^2_\rho \theta + \frac{1}{\rho} \partial_\rho \theta  \right) - \psi_{\text b}''(\tau).
\end{equation}

We employ the Fourier decomposition
\begin{equation}
   \theta(\tau,\rho) = \sum_{n=1}^{\infty} c_n(\tau) J_{0}(j_n \rho),
\end{equation}
where \(J_\nu\) is the Bessel function of the first kind and \(j_n\) is the \(n\)-th root of \(J_{0}\).  
Then \eqref{bessel} yields a system of ODEs
\begin{equation}\label{eq:cn_ode}
   c_n ''(\tau) + (j_n k)^2 c_n (\tau) = - \frac{2}{J^2_{1}(j_n)} \psi_{\text b}''(\tau) \int_0^1 \rho J_{0}(j_n \rho)\,\d \rho,
\end{equation}
where
\begin{equation}
  k  = \sqrt{\frac{1}{4 \sigma_L}}.
\end{equation}
Our initial conditions from \eqref{eq:IBVP_psi} yield initial conditions
\begin{equation}\label{eq:cn_ic}
\begin{aligned}
   &c_n(0) = \psi_{\text b}(0) \frac{2}{J^2_{1}(j_n)}\int_0^1 \rho J_{0}(j_n \rho)\,\d \rho = 0, \\ 
   &c'_n(0) = -\psi_{\text b}'(0) \frac{2}{J^2_{1}(j_n)}\int_0^1 \rho J_{0}(j_n \rho)\,\d \rho = 0.
\end{aligned}
\end{equation}
After one solves the linear system of ordinary differential equations (ODE) (\ref{eq:cn_ode}, \ref{eq:cn_ic})
and performs the inverse substitution for $\psi$ and $\sigma$, one obtains
   \begin{equation}\label{eq:psi_sol}
   \begin{aligned}
      \psi(\sigma, \tau) &= \sum_{n=1}^\infty c_n(\tau) J_0(j_n \sqrt{\sigma / \sigma_L}) + \psi_{\text b}(\tau)
   \end{aligned}
\end{equation}
Similarly, the expression for \(\phi\) is obtained:
\begin{equation}\label{eq:phi_sol}
    \phi(\sigma, \tau) = \frac{1}{2 \sigma_L}\left( \frac{\sigma_L}{\sigma} \right)^{1 / 2} \sum_{n=1}^{\infty} j_n d_n(\tau) J_{1} (j_n \sqrt{\sigma / \sigma_L}),
\end{equation}
where
\begin{equation}
    d_n(\tau) = \int_0^\tau c_n(\xi)\,\d \xi.
\end{equation}
We take the limit as $\sigma \to 0$ in \eqref{eq:psi_sol} and utilise the asymptotic expansion
\begin{equation}
    J_\alpha(z) \sim \frac{1}{\Gamma(\alpha+1)} \left(\frac{z}{2}\right)^\alpha\text{ as } z \to 0.
\end{equation}
From this, we obtain
\begin{equation}\label{discrShorelineeq}
   \psi_{\text{sh}} (\tau)= \psi(0, \tau) =  \sum_{n=1}^\infty c_n(\tau)  + \psi_{\text b}(\tau).
\end{equation}
Once the explicit formula for $c_n(\tau)$ is found, this can be rewritten to obtain what we call the shoreline equation, with $c_n(\tau)$ defined in terms of new constants $a_n$ and $b_n$:
\begin{equation}\label{eq:shoreline}
\begin{aligned}
   \psi_{\text{sh}}(\tau) & =  \sum_{n=1}^\infty
   \frac{b_n }{\sqrt{a_n}} \int_0^ \tau \sin \left( \sqrt{a_n} (\tau- \xi) \right) \psi_{\text{b}}''(\xi) \,\d \xi
   + \psi_{\text b}(\tau)  \\ 
   & = \sum_{n=1}^\infty
   \frac{b_n }{\sqrt{a_n}} 
   \sin \left( \sqrt{a_n} \tau \right) * \psi_{\text{b}}''(\tau)
   + \psi_{\text b}(\tau),
\end{aligned}
\end{equation}
where \(*\) denotes convolution
and 
\begin{equation}
\begin{aligned}
   &a_n =  \frac{j_n^2}{4  \sigma_L}, \\ 
   &b_n =
   -\frac{2}{j_{n} J_{1}\left(j_{n}\right)}.
\end{aligned}
\end{equation}
Equation \eqref{eq:shoreline} solves the forward problem. It is left to note that \(\psi(0,\tau)\) can be easily transformed to the run-up function from the CGT, which at the shoreline reads as
\begin{equation}\label{eq:CG_at_shore}
   {\tau = t+\partial_t R (t)}
   ,\quad
   {\psi(0,\tau) = R (t) + \tfrac{1}{2}(\partial_t R (t))^2}.
\end{equation}
We now proceed to the inverse problem. We take the Laplace transform of \eqref{eq:shoreline} and, by interchanging integration and summation, obtain
\begin{equation}\label{eq:laplace_transform}
\begin{aligned}
    \L \psi_{\text{sh}}(\tau) &= \sum_{n=1}^\infty
    \frac{b_n }{\sqrt{a_n}} 
    \L \sin \left( \sqrt{a_n} \tau \right) \L \psi_{\text{b}}''(\tau) + \L \psi_{\text{b}}(\tau) \\ 
    & = \sum_{n=1}^\infty
    \frac{b_n }{\sqrt{a_n}} 
    \L \sin \left( \sqrt{a_n} \tau \right) s^2 \L \psi_{\text{b}}(\tau) + \L \psi_{\text{b}}(\tau) \\ 
    & = \L \psi_{\text{b}}(\tau) \left(  {s^2}\sum_{n=1}^\infty
    \frac{b_n }{\sqrt{a_n}} 
    \L \sin \left( \sqrt{a_n} \tau \right)  +1 \right) \\ 
    & = \L \psi_{\text{b}}(\tau) \left(  {s^2}\sum_{n=1}^\infty
    \frac{b_n }{a_n + s^2}  +1 \right).
\end{aligned}
\end{equation}
We then regroup the terms and take the inverse Laplace transform of \eqref{eq:laplace_transform}:
\begin{equation}\label{eq:inv_lap}
    \psi_{\text{b}}(\tau) = \L^{-1}\left[ \L \psi_{\text{sh}} \left(  {s^2}\sum_{n=1}^\infty
    \frac{b_n}{a_n + s^2}  +1 \right)^{-1} \right].
\end{equation}
For justification of swapping the integral and summation, see Appendix \ref{subsec:remarks}. It is now left to recover \( u(L,t), \eta (L, t)\) from \(\psi_{\text{b}}(\tau)\). For that, we first compute the solution \(\psi(\sigma, \tau), \phi(\sigma, \tau)\) on the entire domain in the hodograph plane from (\ref{eq:psi_sol}, \ref{eq:phi_sol}) and then perform the inverse CGT. We explain the numerical implementation of this process in Subsection \ref{subsec:num}.

We summarise this section with the complete solution to the inverse problem.
Recall that we want to recover boundary conditions \(u(L, t), \eta(L, t)\) from the observations at the shoreline \(R(t)\). We suggest the following algorithm:
\begin{enumerate}[1.]
    \item Convert \(R(t)\) to \(\psi_{\text{sh}}(\tau)=\psi(0,\tau)\) using CGT at the shore \eqref{eq:CG_at_shore};
    \item Find \(\psi_{\text{b}}(\tau) = \psi(\sigma_L, \tau)\) from \eqref{eq:inv_lap};
    \item Compute \(\psi(\sigma, \tau), \phi(\sigma, \tau)\) from \(\psi_{\text{b}}(\tau)\) using (\ref{eq:psi_sol}, \ref{eq:phi_sol});
    \item Perform the inverse CGT to find \(u(L, t), \eta(L, t)\). 
\end{enumerate}

\subsection{Numerical Verification}\label{subsec:num}
\begin{figure}
    \begin{subfigure}{0.45\textwidth}
        \centering
            \includegraphics[width=\linewidth]{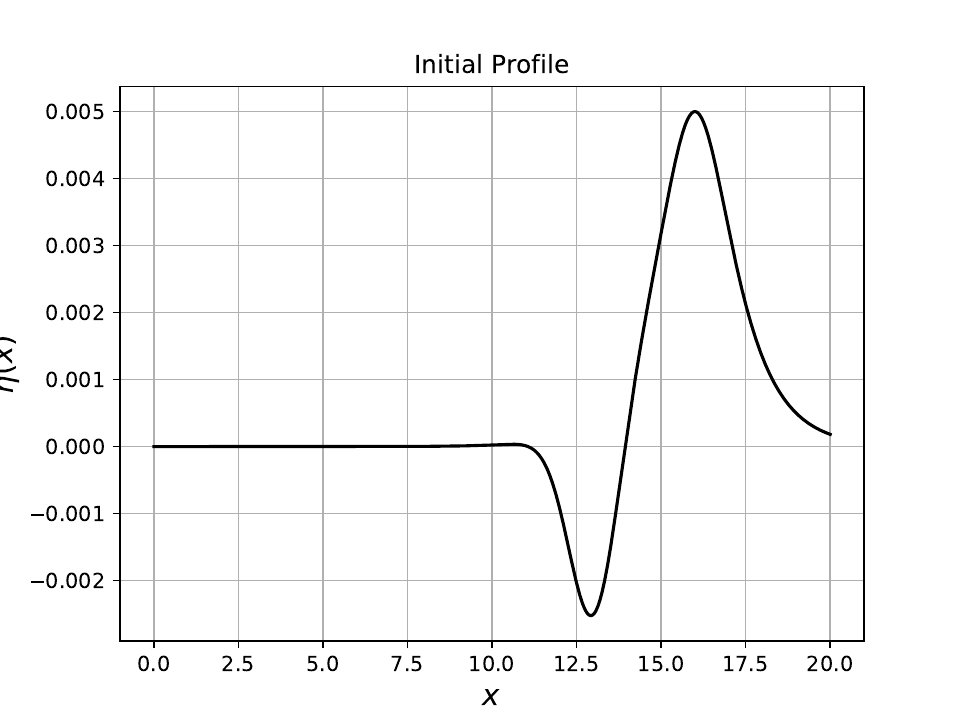}
            \caption{$\eta(x, 0)  = 0.005 \sech(x - 16)- 0.003 e^{-(x-13)^2},$}
            \label{fig:initial_n_wave}
    \end{subfigure}%
        \hfill
    \begin{subfigure}{0.45\textwidth}
            \centering
            \includegraphics[width=\linewidth]{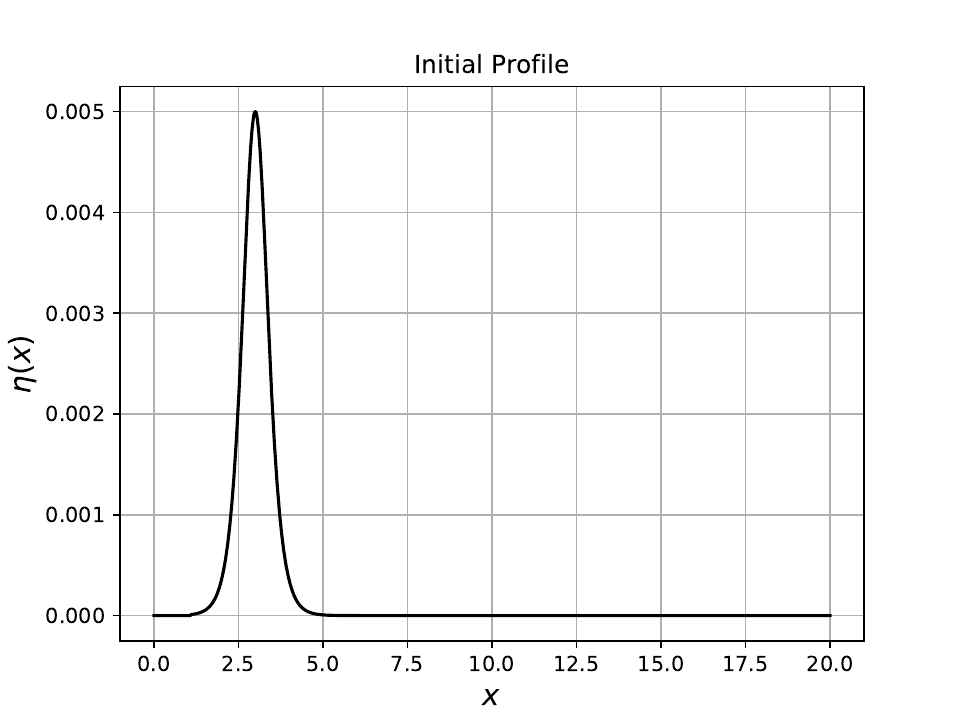}
            \caption{$\eta(x, 0)  = 0.005 \cosh^{-2}(2x-6),$}
            \label{fig:initial_soliton}
    \end{subfigure}%
    \\
    \begin{subfigure}{0.45\textwidth}
        \centering
            \includegraphics[width=\linewidth]{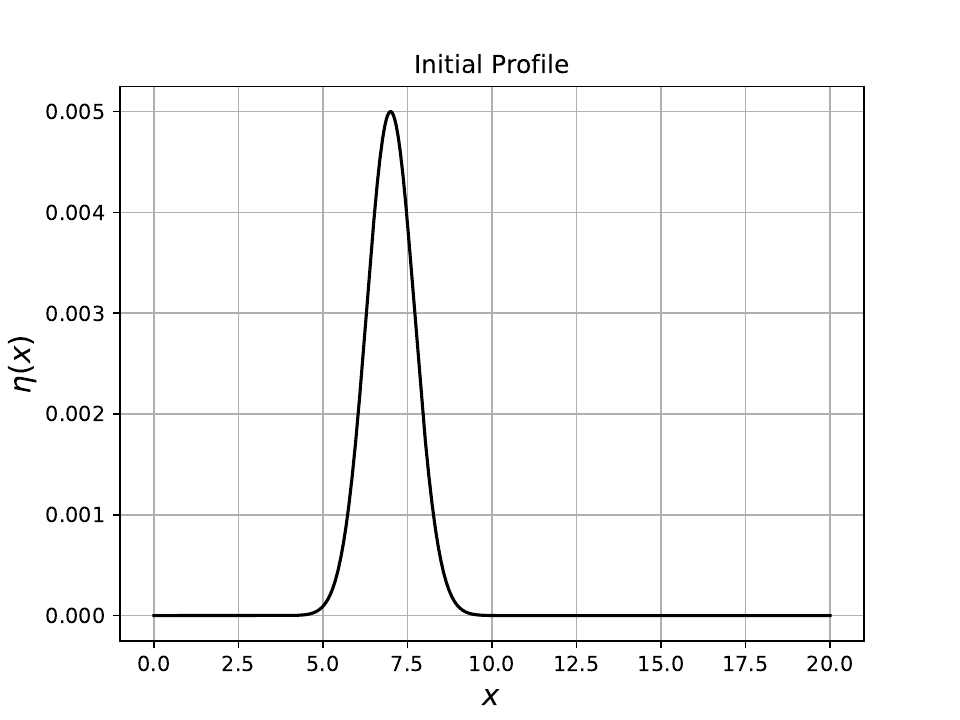}
            \caption{$\eta(x, 0)  = 0.005 e^{-(x-7)^2},$}
            \label{fig:initial_gauss}
    \end{subfigure}%
        \hfill
    \begin{subfigure}{0.45\textwidth}
            \centering
            \includegraphics[width=\linewidth]{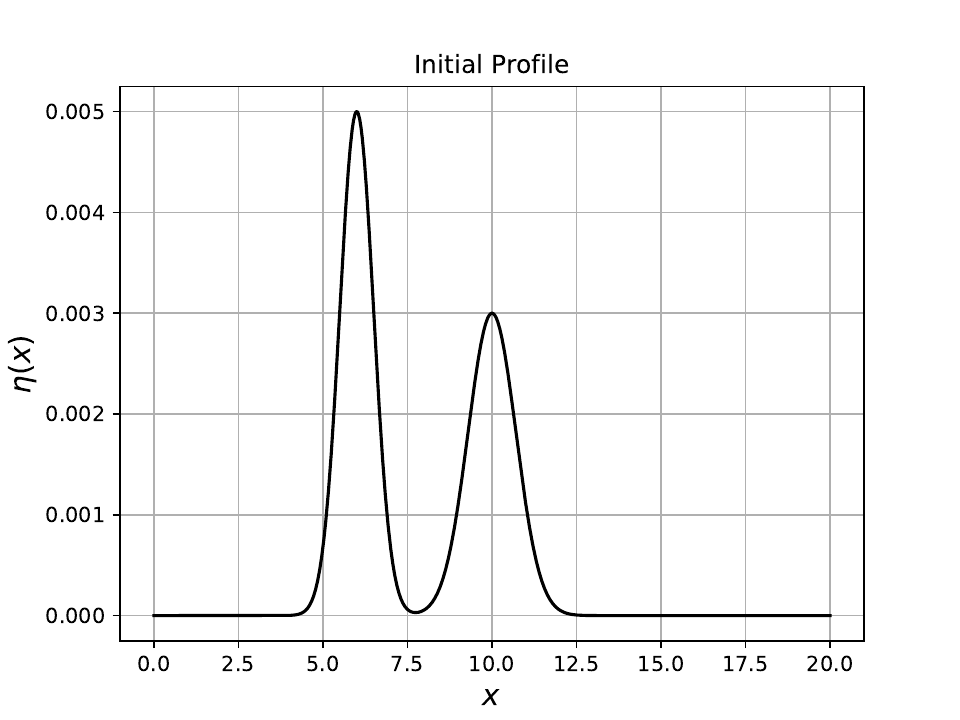}
            \caption{$\eta(x, 0)  = 0.005  e^{-2  (x-6)^2 }  + 0.003 e^{-(x-10) ^ 2}$.}
            \label{fig:initial_2_soliton}
    \end{subfigure}%
    \caption{Initial displacements as defined in \eqref{eq:ics_eta0}.}\label{fig:eta0}
\end{figure}

To confirm the above formulae, we verify our proposed algorithm numerically. For convenience, we set the boundary of the finite sloped interval to be $L=1$. Additionally, note that all verification was performed in dimensionless units. To dimensionalise, a scaling factor can be added and is given in \eqref{eq:dimensions}.
We manufacture a solution to the inverse problem. A summary of this process can be seen in the diagram below.
\[\begin{tikzcd}
	{\eta_0} & {\psi_0} & {\psi_{\text{sh}}} & {\psi_{\text{b}}} & {\psi,\phi} & {\eta_{\text{b}}, u_{\text{b}}} \\
	&&& {\psi_{\text{b}}}
	\arrow[from=1-1, to=1-2]
	\arrow[from=1-2, to=1-3]
	\arrow[curve={height=12pt}, from=1-2, to=2-4]
	\arrow[from=1-3, to=1-4]
	\arrow[from=1-4, to=1-5]
	\arrow[from=1-5, to=1-6]
	\arrow[shift left, no head, from=2-4, to=1-4]
	\arrow[no head, from=2-4, to=1-4]
\end{tikzcd}\]
Consider an infinitely-sloping plane beach with an initial value problem on the \( x \) semi-axis. We choose the initial conditions to be
\begin{equation}\label{eq:ics_eta0}
\begin{aligned}
    \eta(x, 0) & = 0.005 \sech(x - 16)- 0.003\exp(-(x-13)^2) ,\\
    \eta(x, 0) & = 0.005 \cosh^{-2}(2x-6), \\ 
    \eta(x, 0) & = 0.005 \exp(-(x-7)^2), \\ 
    \eta(x, 0) & = 0.005  \exp(-2  (x-6)^2 )  + 0.003 \exp(-(x-10) ^ 2)
\end{aligned}
\end{equation}
with zero initial velocity
\begin{equation}\label{eq:u_0}
    u(x, 0) = 0
\end{equation}
to simulate an earthquake as the initial disturbance. Additionally, it provides greater simplicity for the problem as this case does not require data projection.
Initial profiles can be seen in Figure \ref{fig:eta0}.

\begin{figure}
    \begin{subfigure}{0.45\textwidth}
        \centering
            \includegraphics[width=\linewidth]{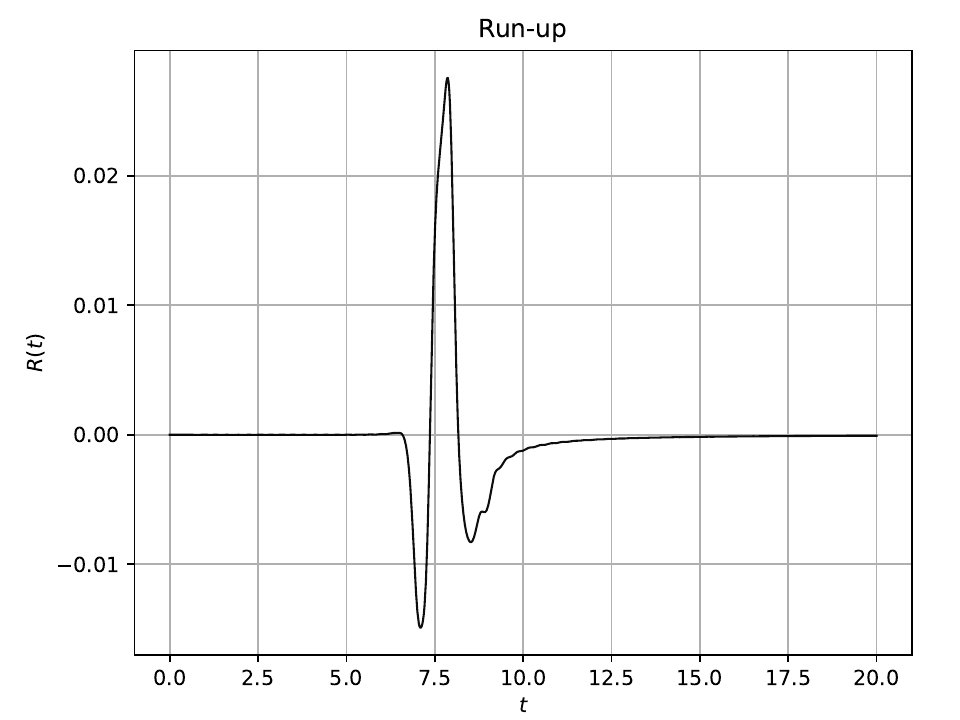}
            \caption{$\eta(x, 0)  = 0.005 \sech(x - 16)- 0.003 e^{-(x-13)^2},$}
            \label{fig:runup_n_wave}
    \end{subfigure}%
        \hfill
    \begin{subfigure}{0.45\textwidth}
            \centering
            \includegraphics[width=\linewidth]{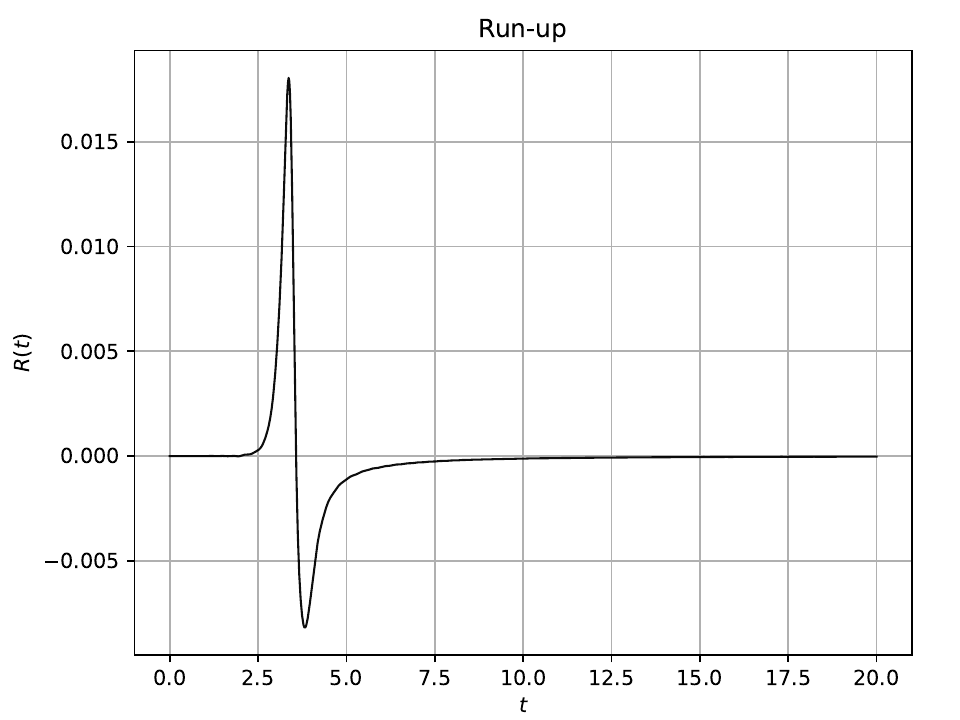}
            \caption{$\eta(x, 0)  = 0.005 \cosh^{-2}(2x-6),$}
            \label{fig:runup_soliton}
    \end{subfigure}%
    \\
    \begin{subfigure}{0.45\textwidth}
        \centering
            \includegraphics[width=\linewidth]{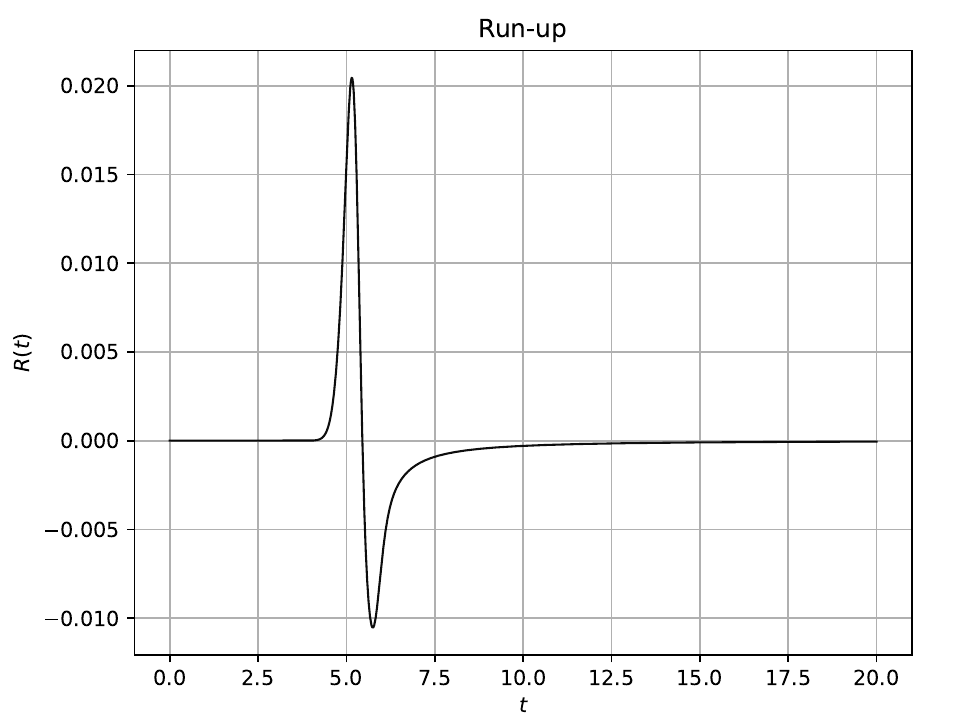}
            \caption{$\eta(x, 0)  = 0.005 e^{-(x-7)^2},$}
            \label{fig:runup_gauss}
    \end{subfigure}%
        \hfill
    \begin{subfigure}{0.45\textwidth}
            \centering
            \includegraphics[width=\linewidth]{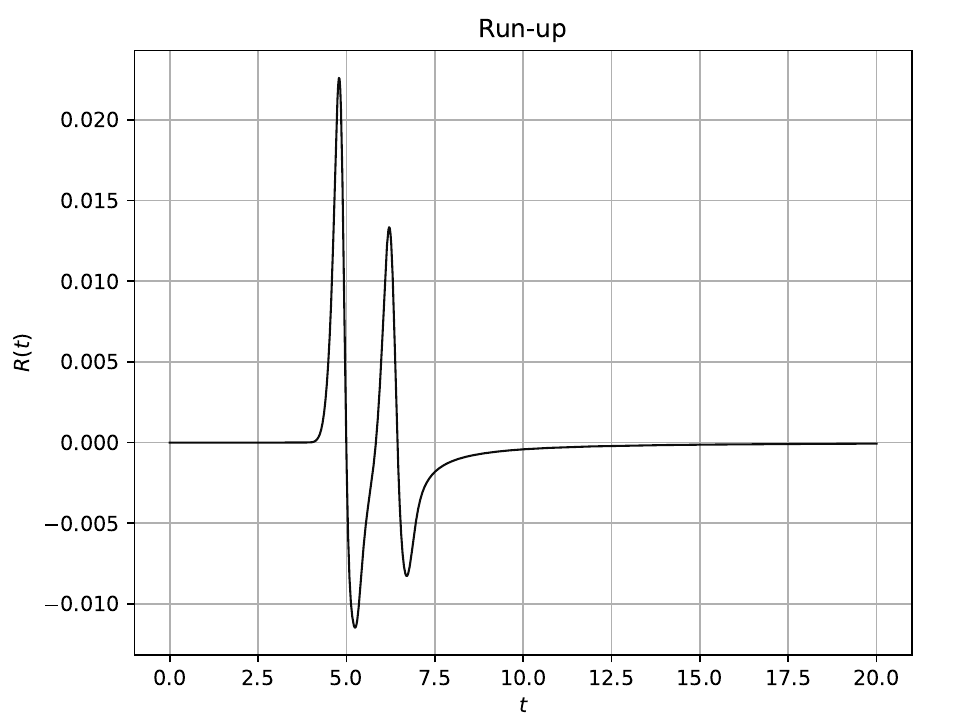}
            \caption{$\eta(x, 0)  = 0.005  e^{-2  (x-6)^2 }  + 0.003 e^{-(x-10) ^ 2}$.}
            \label{fig:runup_2_soliton}
    \end{subfigure}%
    \caption{Run-ups corresponding to IC in \eqref{eq:ics_eta0}.}\label{fig:runups}
\end{figure}

\begin{figure}
    \begin{subfigure}{0.45\textwidth}
        \centering
            \includegraphics[width=\linewidth]{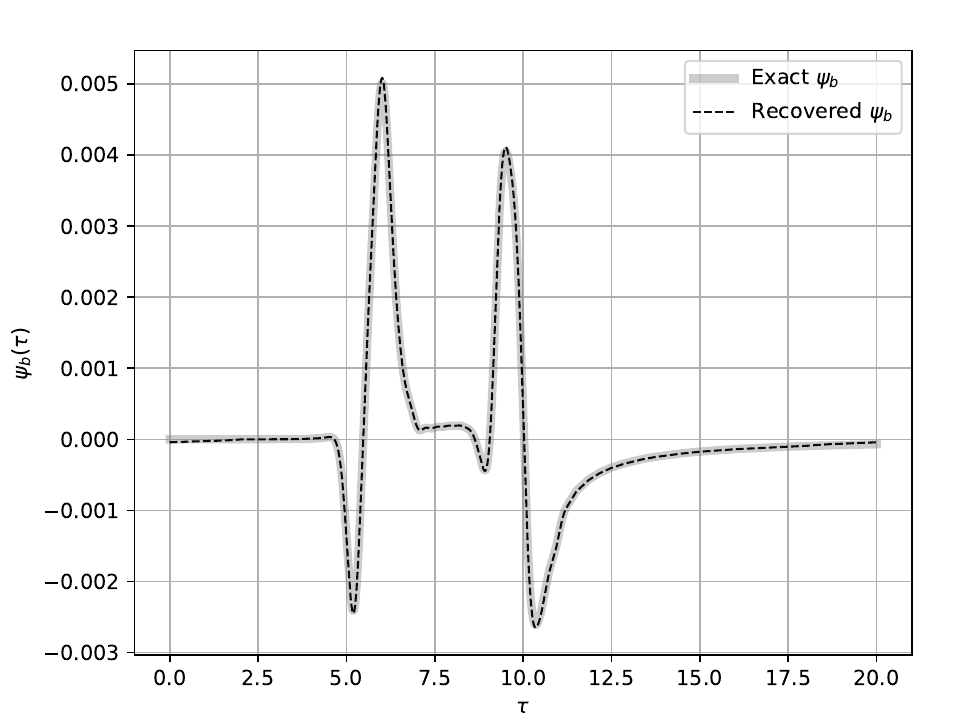}
            \caption{$\eta(x, 0)  = 0.005 \sech(x - 16)- 0.003 e^{-(x-13)^2},$}
            \label{fig:psib_sech_gauss}
    \end{subfigure}%
        \hfill
    \begin{subfigure}{0.45\textwidth}
            \centering
            \includegraphics[width=\linewidth]{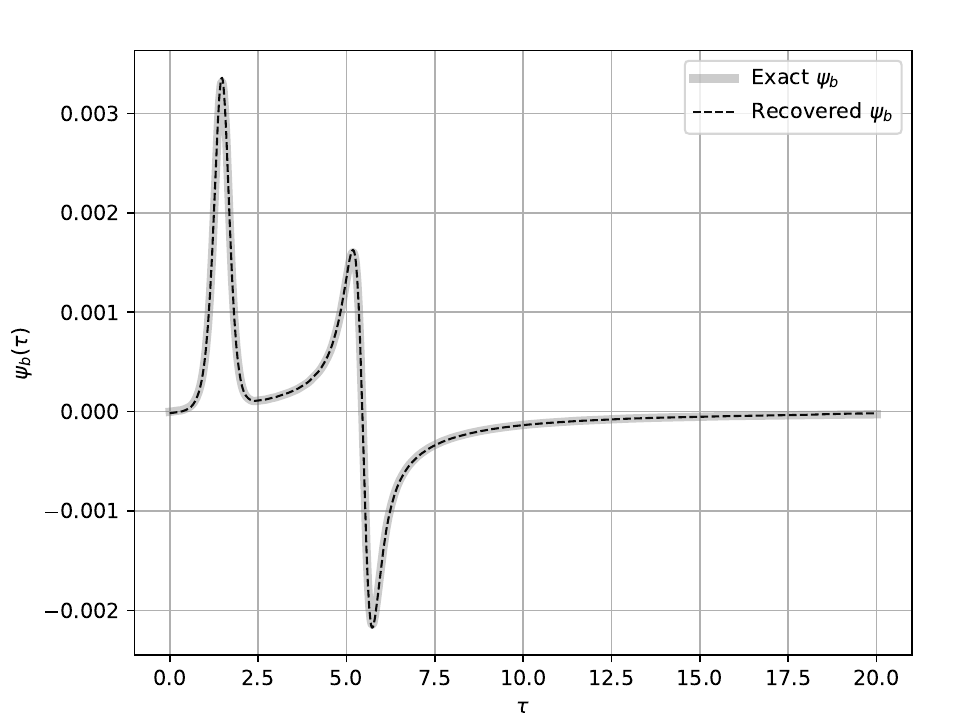}
            \caption{$\eta(x, 0)  = 0.005 \cosh^{-2}(2x-6),$}
            \label{fig:psib_soliton}
    \end{subfigure}%
    \\
    \begin{subfigure}{0.45\textwidth}
        \centering
            \includegraphics[width=\linewidth]{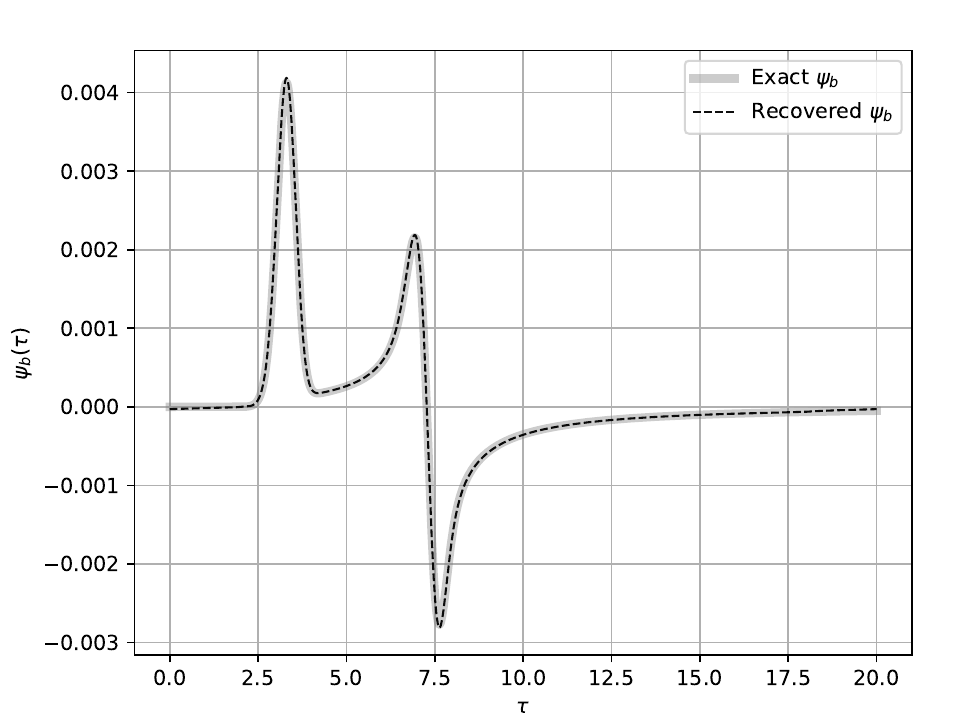}
            \caption{$\eta(x, 0)  = 0.005 e^{-(x-7)^2},$}
            \label{fig:psib_gauss}
    \end{subfigure}%
        \hfill
    \begin{subfigure}{0.45\textwidth}
            \centering
            \includegraphics[width=\linewidth]{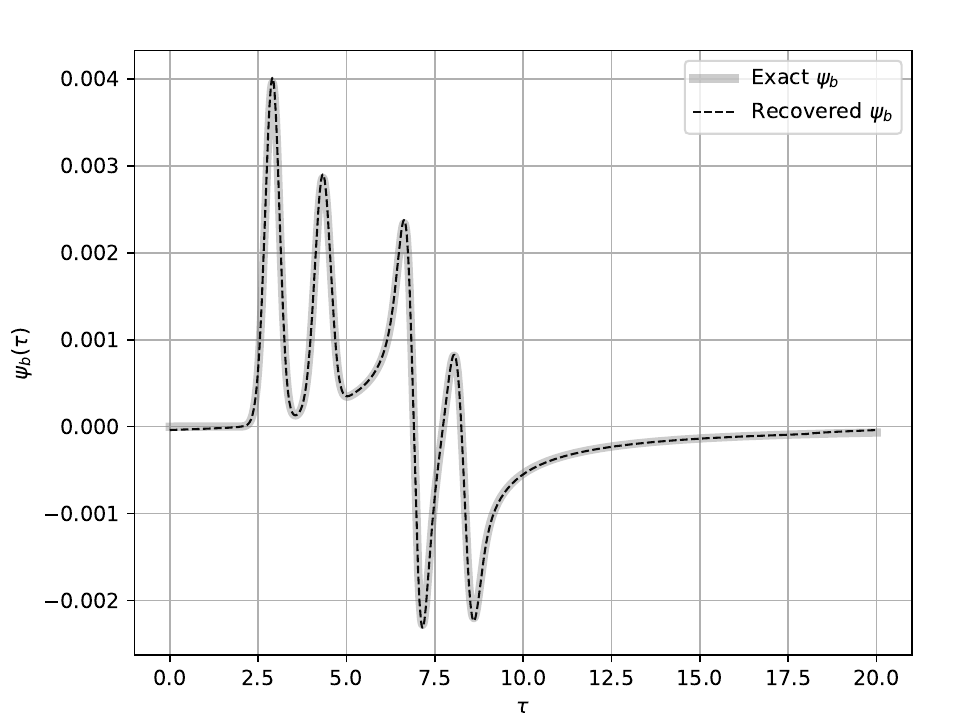}
            \caption{$\eta(x, 0)  = 0.005  e^{-2  (x-6)^2 }  + 0.003 e^{-(x-10) ^ 2}$.}
            \label{fig:psib_2_soliton}
    \end{subfigure}%
    \caption{Comparison of \(\psi_{\text{b}}(\tau)\) and \(\psi_{\text{b}}^{\text{e}}(\tau)\) corresponding to IC in \eqref{eq:ics_eta0}.}\label{fig:psi_compared}
\end{figure}

\begin{figure}
    \begin{subfigure}{0.45\textwidth}
        \centering
            \includegraphics[width=\linewidth]{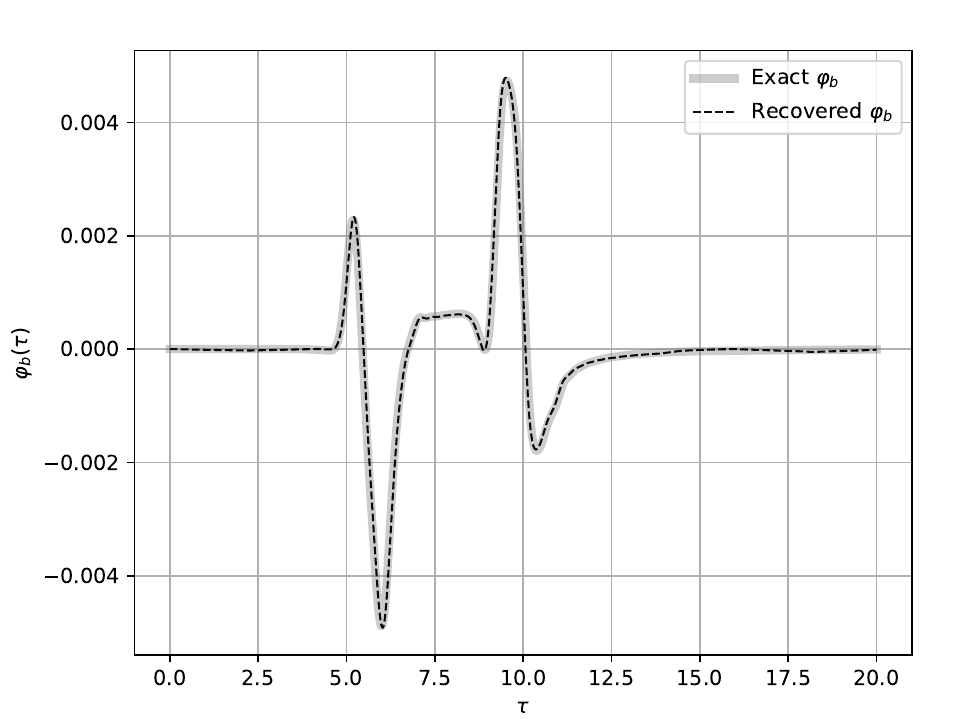}
            \caption{$\eta(x, 0)  = 0.005 \sech(x - 16)- 0.003 e^{-(x-13)^2},$}
            \label{fig:phib_n_wave}
    \end{subfigure}%
        \hfill
    \begin{subfigure}{0.45\textwidth}
            \centering
            \includegraphics[width=\linewidth]{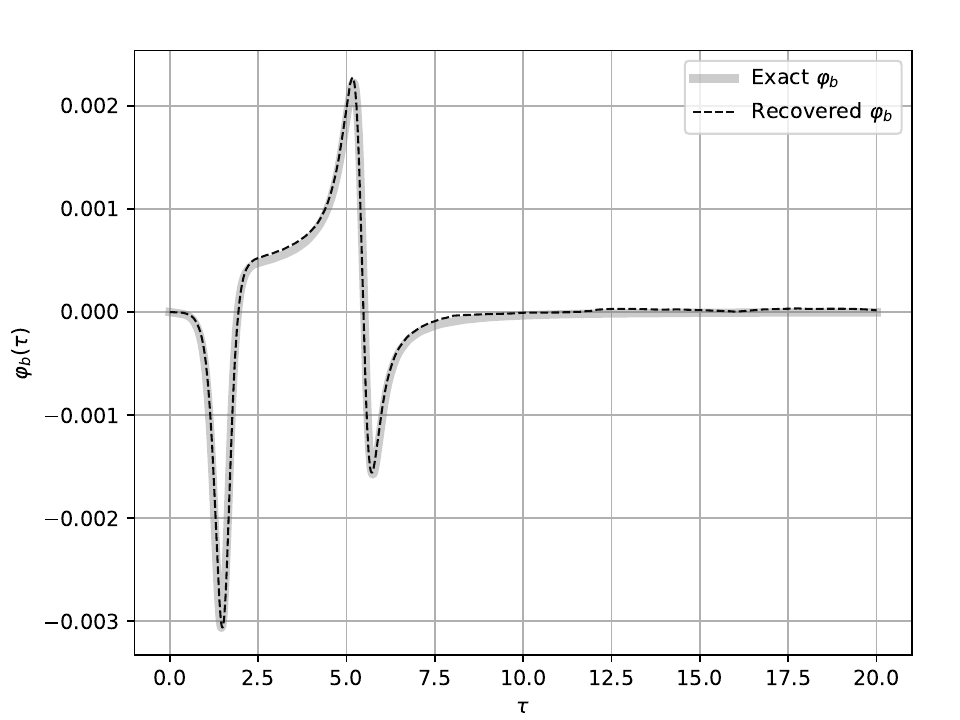}
            \caption{$\eta(x, 0)  = 0.005 \cosh^{-2}(2x-6),$}
            \label{fig:phib_soliton}
    \end{subfigure}%
    \\
    \begin{subfigure}{0.45\textwidth}
        \centering
            \includegraphics[width=\linewidth]{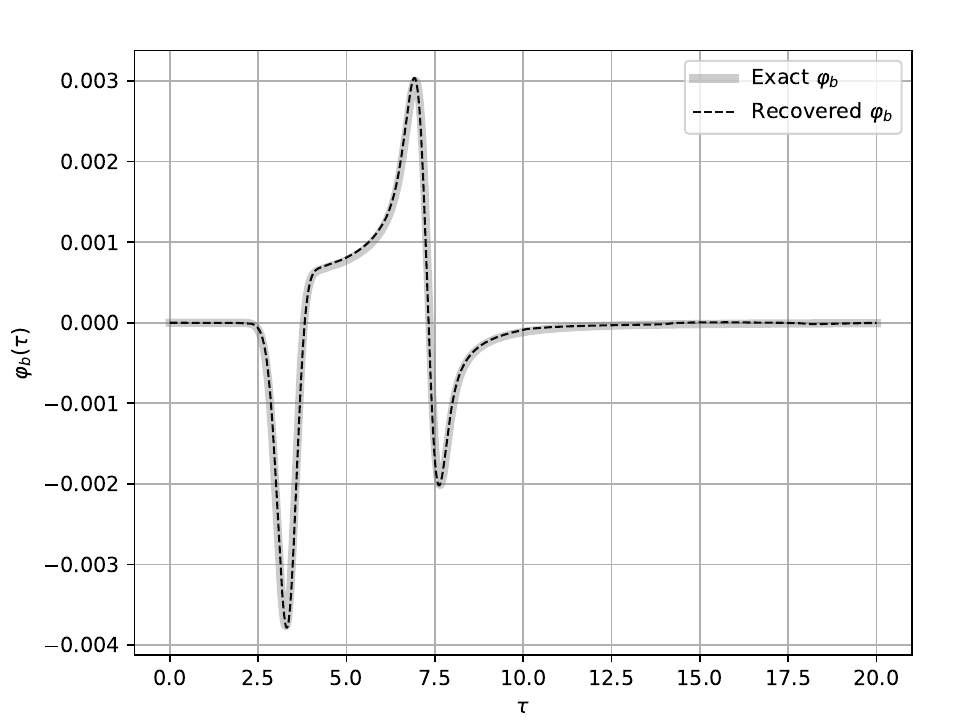}
            \caption{$\eta(x, 0)  = 0.005 e^{-(x-7)^2},$}
            \label{fig:phib_gauss}
    \end{subfigure}%
        \hfill
    \begin{subfigure}{0.45\textwidth}
            \centering
            \includegraphics[width=\linewidth]{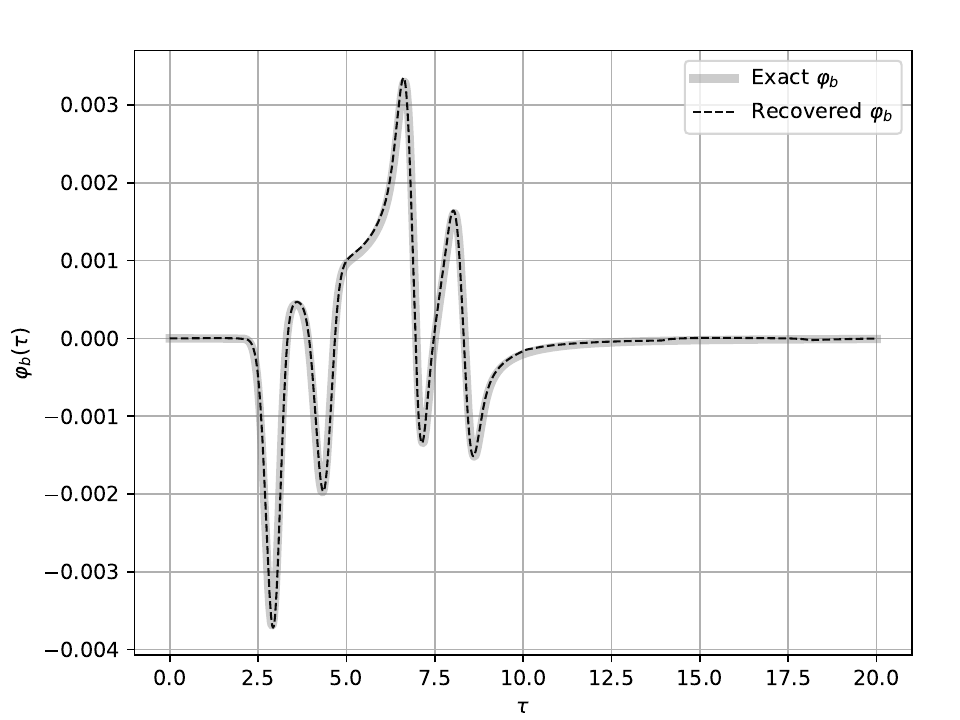}
            \caption{$\eta(x, 0)  = 0.005  e^{-2  (x-6)^2 }  + 0.003 e^{-(x-10) ^ 2}$.}
            \label{fig:phib_2_soliton}
    \end{subfigure}%
    \caption{Comparison of \(\phi_{\text{b}}(\tau)\) and \(\phi_{\text{b}}^{\text{e}}(\tau)\) corresponding to IC in \eqref{eq:ics_eta0}.}\label{fig:phi_compared}
\end{figure}

We solve the forward problem (IVP (\ref{eq:SWE}, \ref{eq:ics_eta0}, \ref{eq:u_0}) on \(x>0\), \(t> 0\)) using techniques of \cite{Rybkin21} to find the run-up \(R(t)\), and while simulating the run-up we also record \(\psi_{\text{b}}^{\text{e}}(\tau) = \psi(1, \tau)\) and \(\phi_{\text{b}}^{\text{e}}(\tau) = \phi(1, \tau)\). Superscript e stands for exact.
Corresponding run-ups can be seen in Figure \ref{fig:runups}.
Next, our algorithm outlined above is applied to \(R(t)\).
First, we compare the found \(\psi_{\text{b}}(\tau)\) and \(\phi_{\text{b}}(\tau)\) to \(\psi_{\text{b}}^{\text{e}}(\tau)\) and \(\phi_{\text{b}}^{\text{e}}(\tau)\) respectively. Excellent agreement can be seen in Figures \ref{fig:psi_compared} and \ref{fig:phi_compared}.

\begin{figure}
    \begin{subfigure}{0.45\textwidth}
        \centering
            \includegraphics[width=\linewidth]{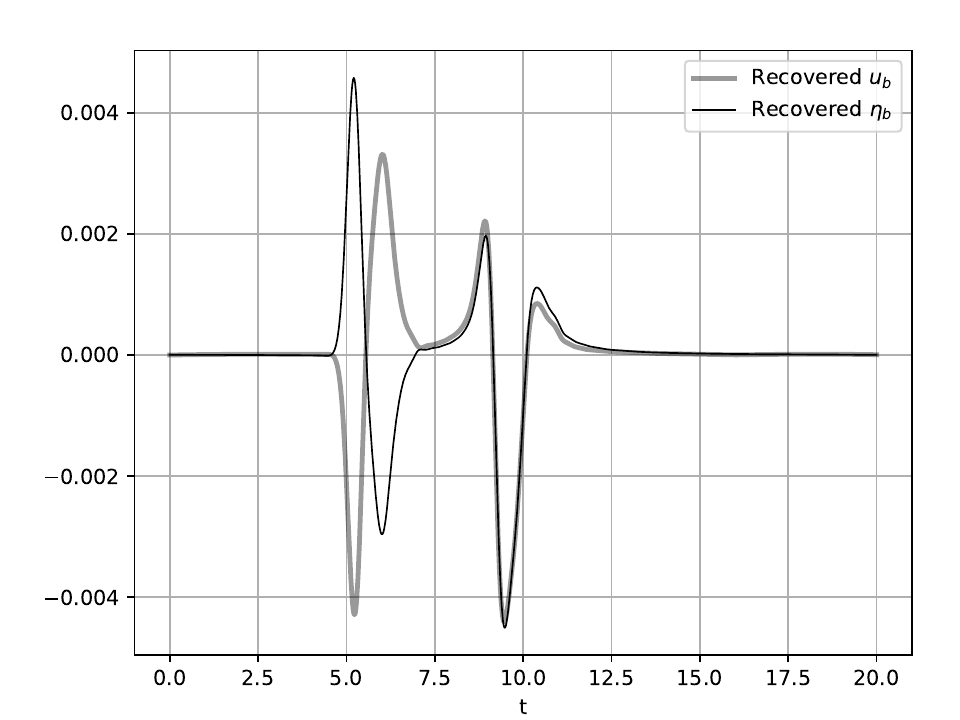}
            \caption{$\eta(x, 0)  = 0.005 \sech(x - 16)- 0.003 e^{-(x-13)^2},$}
            \label{fig:bc_sech_gauss}
    \end{subfigure}%
        \hfill
    \begin{subfigure}{0.45\textwidth}
            \centering
            \includegraphics[width=\linewidth]{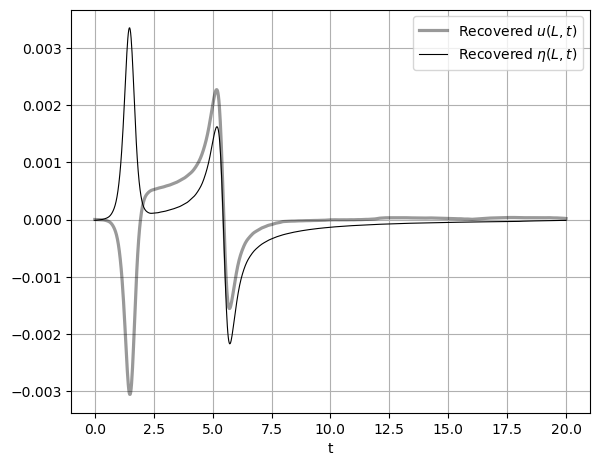}
            \caption{$\eta(x, 0)  = 0.005 \cosh^{-2}(2x-6),$}
            \label{fig:bc_soliton}
    \end{subfigure}%
    \\
    \begin{subfigure}{0.45\textwidth}
        \centering
            \includegraphics[width=\linewidth]{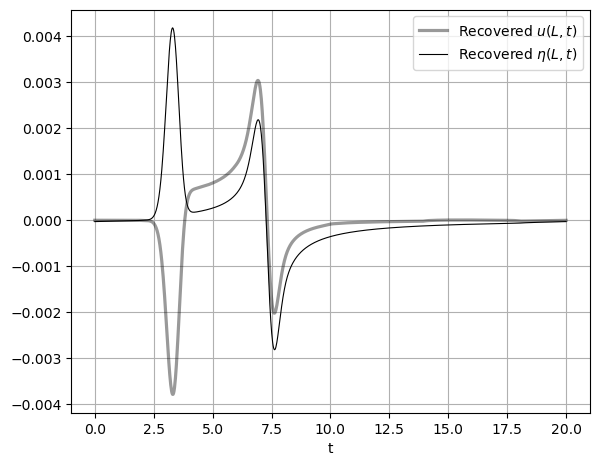}
            \caption{$\eta(x, 0)  = 0.005 e^{-(x-7)^2},$}
            \label{fig:bc_gauss}
    \end{subfigure}%
        \hfill
    \begin{subfigure}{0.45\textwidth}
            \centering
            \includegraphics[width=\linewidth]{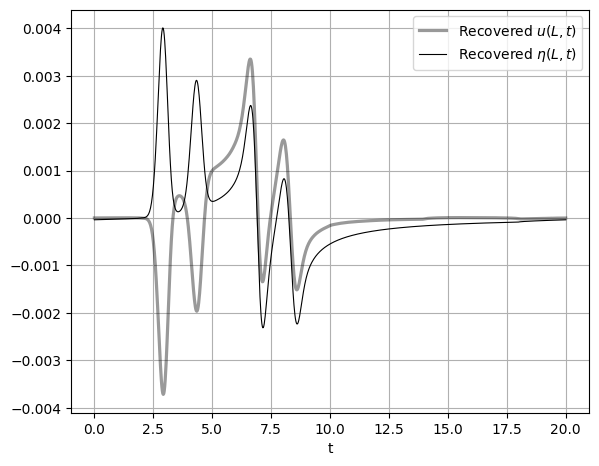}
            \caption{$\eta(x, 0)  = 0.005  e^{-2  (x-6)^2 }  + 0.003 e^{-(x-10) ^ 2}$.}
            \label{fig:bc_two_bumps}
    \end{subfigure}%
    \caption{Recovered boundary conditions corresponding to IC in \eqref{eq:ics_eta0}.}\label{fig:bcs_recovered}
\end{figure}

We proceed with our algorithm. Once \(\psi(\sigma, \tau)\) and \(\phi(\sigma, \tau)\) are found, it is left to perform the inverse CGT. As noted, the physical initial conditions \(u(L, t), \eta(L, t)\) are specified on a parametric curve in the \((\sigma, \tau)\) plane, \(\Gamma = (\sigma(\tau), \tau)\). The CGT on \(\Gamma\) reads as 
\begin{equation}
    \sigma = L + \psi - \tfrac{1}{2}\phi^2, \quad \tau = t - \phi.
\end{equation}
So we compute \(\psi, \phi\) on a grid and then for each value of \(\tau\) we find \(\sigma\) so that \( \abs{\sigma - L - \psi + \phi^2 / 2} \) is minimised in order to find \(\sigma(\tau)\). Finally, we recover
\begin{equation}
     \eta(L,t) = \restr{\psi}{\Gamma} - \tfrac{1}{2}\restr{\phi^{2}}{\Gamma}, \quad u(L,t) = \restr{\phi}{\Gamma}, \quad t = \tau + \restr{\phi}{\Gamma}. 
\end{equation}
Recovered physical boundary conditions can be seen in Figure \ref{fig:bcs_recovered}.
\subsection{Computational Cost}
We briefly discuss the computational cost of our proposed algorithm.
Given the run-up data $R(t)$ as an $N$-dimensional vector of measurements, we compute $u(L,t)$ and $\eta(L,t)$ using various numerical methods supported by the algorithm described in a previous section. The implementation begins with the Carrier-Greenspan transform in \eqref{CG} which requires $O(N)$ operations to obtain $\psi(0,\tau)$ and $\phi(0,\tau)$. 
Next, we apply the Laplace transform in \eqref{eq:inv_lap} at a cost of $O(N^2)$. Implementing the inverse Laplace transform to obtain $\psi_{\text{b}}(\tau)$ in \eqref{eq:inv_lap} would typically require $O(N^3\log{N})$. However, under the reasonable assumption that the real component of all poles of the Laplace transformed $\psi_{\text{sh}}(\tau)$ are strictly negative, we can instead use the inverse Fast Fourier Transform (iFFT) to significantly improve efficiency as the cost is reduced to $O(N\log{N})$. 
Note that $a_n, b_n, j_n $ are assumed to be precomputed and stored in arrays prior to the summation in \eqref{eq:inv_lap}, so no significant cost is added to the complexity. 
Next, we recover $\psi(\sigma, \tau)$ and $\phi(\sigma, \tau)$ using (\ref{eq:psi_sol}, \ref{eq:phi_sol}) which is the most computationally expensive at $O(N^3)$. This is primarily due to the calculation of the coefficients $c_n$, $d_n$ involved in determining $\psi,\phi$. 
To apply the inverse Carrier-Greenspan transform, the parametric curve $\Gamma$ must be recovered through minimisation, which is done via brute force with a complexity of $O(N^2)$. This recovers $u(L,t)$ and $\eta(L,t)$.  
Thus, the computational cost for the entire algorithm is $O(N^3)$.

\begin{figure} [H]
    \centering
    \includegraphics[width=0.75\linewidth]{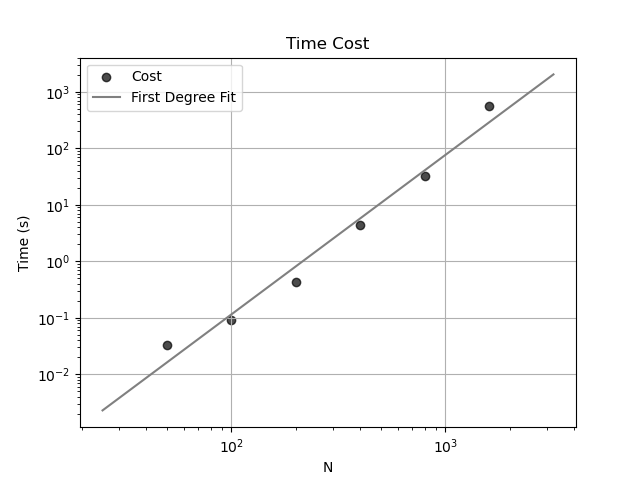}
    \caption{Least square fit of the time with a polynomial of degree 1. Slope is 2.82.}
    \label{fig:time_N}
\end{figure}

It took 76 seconds to run our algorithm with 1500 data points and truncating $n$ at 500 on Google Colab (here \(n\) is the number of terms in the truncated series from \eqref{eq:inv_lap}). We also measured the time as a function of the number of points $N$ on our laptop (Figure \ref{fig:time_N}).

\section{Boussinesq Equation} \label{sec:bouss}
\begin{figure}
    \centering
    \includegraphics[width=0.5\linewidth]{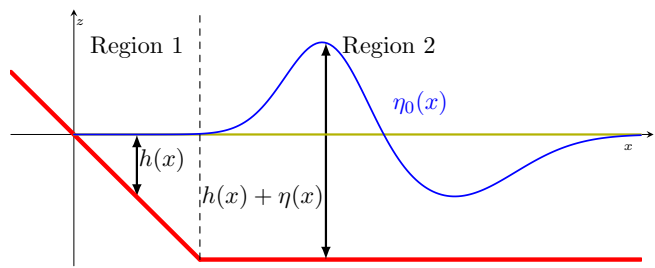}
    \caption{Cross-sectional view of the piece-wise sloping bathymetry. The virtual buoy is located at the intersection of the dashed and blue lines.}
    \label{fig:bath_complex}
\end{figure}

Now we discuss stitching the NSWE model with the Boussinesq equation to examine Region 2 of Figure \ref{fig:bath_complex}. To accomplish the stitching, we will first discuss the process heuristically. Then, with this groundwork, we will describe the actual process and implementation.
Let us first consider the classical (``bad'') Boussinesq equation. It is considered to be ``bad'', from the mathematical perspective, because it is linearly ill-posed; it differs from the ``good'' Boussinesq equation by the sign in front of the fourth derivative. The bad Boussinesq equation is fundamental for modelling non-linear, dispersive water waves \cite{Johnson97}. It is a completely integrable equation and permits soliton solutions. These solutions can be obtained using a variety of methods, including the inverse scattering transform \cite{deift82} and Hirota's bilinear method \cite{hirota73, Johnson97}. On the half-line, the unified transform method can also be used to obtain integral-type solutions \cite{Fokas}. 
The ``good'' Boussinesq equation finds its applications in modelling non-linear vibrations along a string \cite{zakharov73}.
\subsection{Heuristic Understanding} \label{subsec:bouss_heur}
The key idea is to recover data on the boundary between Regions 1 and 2 from the shoreline data using the result of Section \ref{sec:finite_int}, fit this data to a known multi-soliton solution to the Boussinesq equation at \(x = L\), and then recover the multi-soliton initial conditions in Region 2.

\begin{figure}[H]
    \centering
    \includegraphics[width=0.5\linewidth]{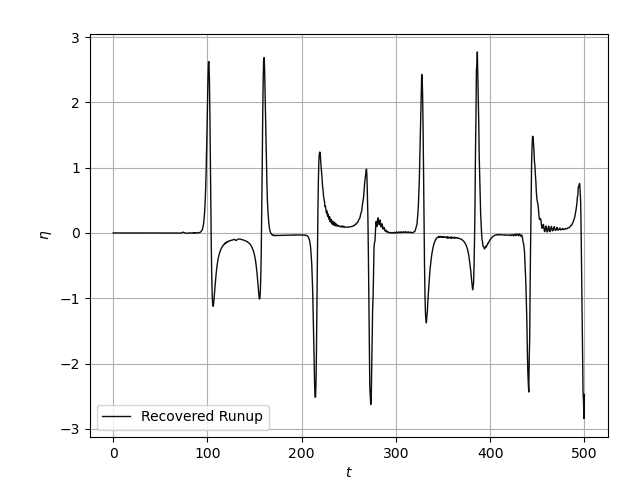}
    \caption{Run-up that correspond to \(\psi_{\text{b}}\) being one soliton. Notice the alternating pattern.}
    \label{fig:Runup_036_010_far_longtime}
\end{figure}

The classical Boussinesq equation is given in dimensionless units as
\begin{equation}\label{eq:boussinesq}
   \partial^2_t \eta = \partial^2_x \eta + 6 \partial_x(\eta^2) + \partial^4_x \eta.
\end{equation}
The natural recourse is to dimensionalise the equation into common units using the substitution:
\begin{equation}
     \wt x = (H_0 / \sqrt{3})x ,\quad \wt t = \sqrt{H_0 / 3g}\, t,\quad \wt \eta = 4H_0 \eta,
\end{equation}
where the tilde signifies the quantity in dimensional units, \(H_0 \) is the characteristic height, and \(g\) is acceleration due to gravity. This ensures that the equations will describe the same physical characteristics. 
Stitching the two equations together requires the Boussinesq equation on the half-line. Solutions and discussions on this were done in \cite{Fokas}. However, for a simpler model, we can assume a solution to the Boussinesq equation on the full line. We want to manufacture a solution to the inverse problem. For a moment, think about the forward problem. We can have a travelling wave (think of a single soliton for now) that passes through the boundary between Regions 1 and 2 in the direction of Region 1. We then read off the data at the boundary and enter it as boundary conditions for the sloping domain (Region 1). Notably, we assume $\psi_{\text{b}} \approx \eta_{\text{b}}$ (equivalently $u_{\text{b}}^2 \ll 1$) as a further simplification (because we want to avoid projections). This assumption is typical: while the speed of propagation of a long tsunami wave over deep water is of order of hundreds metres per second, the flow velocity is of order of tens metres per second; similar arguments have previously been made by \cite{Synolakis87}.
As seen in Figure \ref{fig:bcs_recovered}, indeed \(\eta_{\text{b}}\) is of order \(10^{-3}\)  while \(u_{\text{b}}^2\) is of order \(10^{-6}\). 
Since the boundary condition obtained from the Boussinesq equation does not account for the wave reflected from the shore, this boundary condition (single soliton) leads to the wave ``bouncing'' within the interval, since the boundary value is zero after the wave first enters. An example of this behaviour is given in Figure \ref{fig:Runup_036_010_far_longtime}. 
For the inverse problem, we need shoreline data that corresponds to boundary data obtained from the Boussinesq equation. To generate this, we assume
the initial condition in the Boussinesq region to be the sum of two solitons (multi-soliton solution can be treated in the same manner).
Then we take \(\eta_{\text{b}}\) (under our assumptions same as $\psi_{\text{b}}$) corresponding to this IC. Due to the reflected from the shore wave, we need $L$ large enough that we can distinguish between the wave we set at the boundary and the reflected one. Now, we can either crop the shoreline data \(R(t)\) or boundary data $\psi_{\text{b}}(\tau)$, after recovering it. Since $\psi_{\text{b}}$ does not immediately reflect a wave, it is easier to crop this data (Figure \ref{fig:Recovery_010_036_close}). Note that cropping happens because we need to manufacture appropriate data. When appropriate data is given, no cropping is needed. Using the cropped boundary data, we can extract the parameters for the two solitons using methods such as least squares (Figure \ref{fig:Fitting_010_036_close}), which provides their speeds.
For this part we assume the solution in the Boussinesq region to be a 2-soliton of the form \eqref{2soliton} and then recover the parameters.
Once the solution in the Boussinesq region is obtained, we can recover the initial (at the time of the tsunamigenic event) conditions in the Boussinesq region.

\subsection{Implementation} \label{subsec:bouss_implementation}
\begin{figure}
    \begin{subfigure}{0.45\textwidth}
        \centering
            \includegraphics[width=\linewidth]{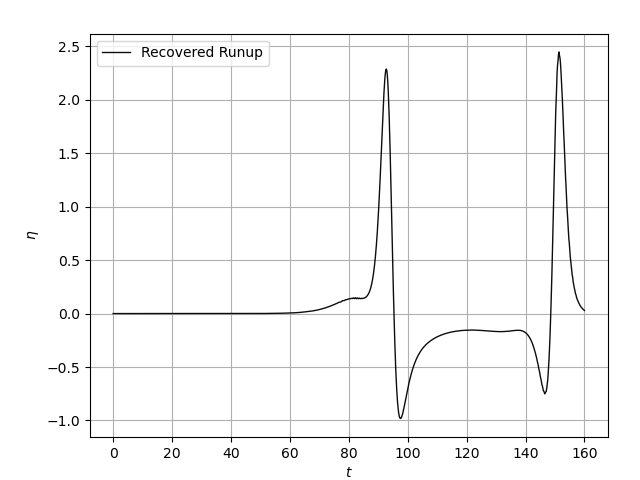}
        \caption{Run-up generated by a two-soliton solution.}
        \label{fig:Runup_010_036_close}
    \end{subfigure}
    \hfill
    \begin{subfigure}{0.45\textwidth}
        \centering
            \includegraphics[width=\linewidth]{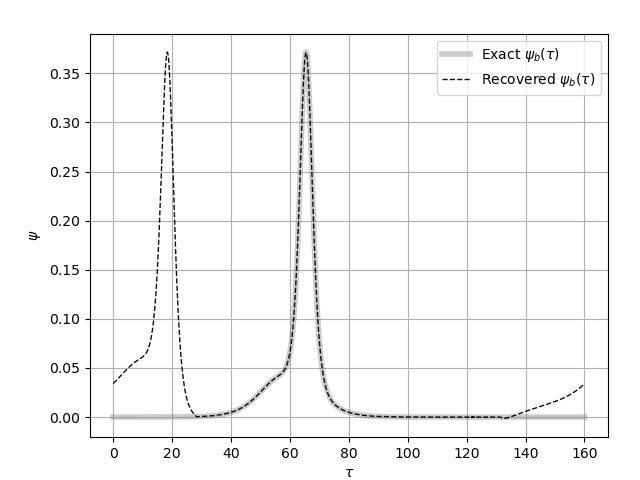}
        \caption{Recovered $\psi_{\text{b}}$ from the run-up.}
        \label{fig:Recovery_010_036_close}
    \end{subfigure}
    \\
    \centering
    \begin{subfigure}{0.45\textwidth}
            \includegraphics[width=\linewidth]{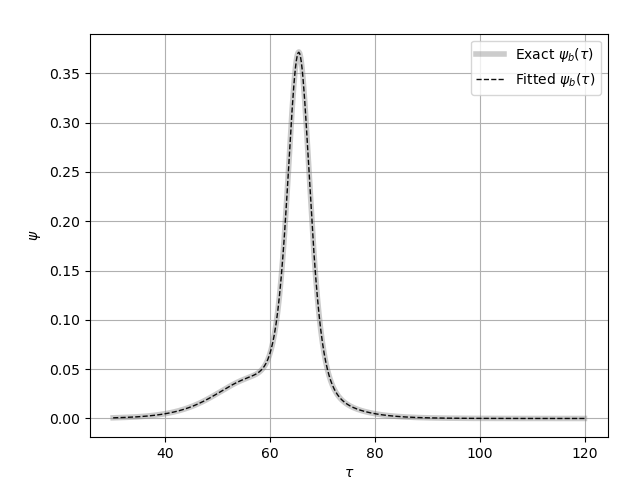}
        \caption{Recovered $\psi_{\text{b}}$ compared with the original Psi boundary: Cropping was at $t = 30$ and $t = 120$.}
        \label{fig:Fitting_010_036_close}
    \end{subfigure}
    \caption{Numerical experiment with soliton parameters $q_1 = 0.1$, $q_2 = \sqrt{0.1}$, $t_1 = 60$, and $t_2 = 65$.}\label{fig:stitch_exp_1}
\end{figure}

\begin{figure}
    \begin{subfigure}{0.45\textwidth}
        \centering
            \includegraphics[width=\linewidth]{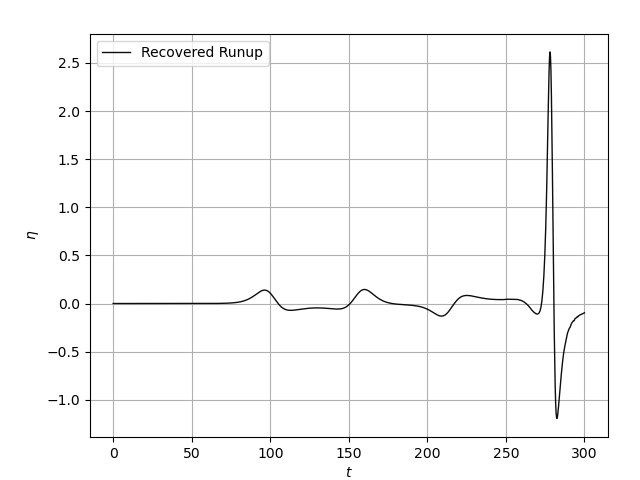}
        \caption{Run-up generated by a two-soliton solution.}
        \label{fig:Runup_010_036_far}
    \end{subfigure}
    \hfill
    \begin{subfigure}{0.45\textwidth}
        \centering
            \includegraphics[width=\linewidth]{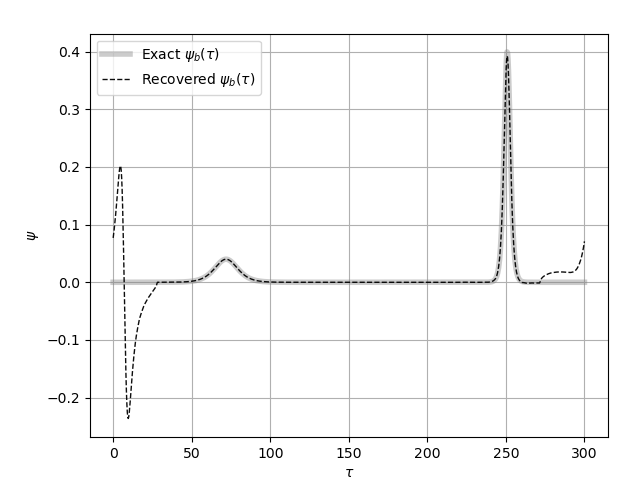}
        \caption{Recovered $\psi_{\text{b}}$ from the run-up.}
        \label{fig:Recovery_010_036_far}
    \end{subfigure}
    \\
    \centering
    \begin{subfigure}{0.45\textwidth}
            \includegraphics[width=\linewidth]{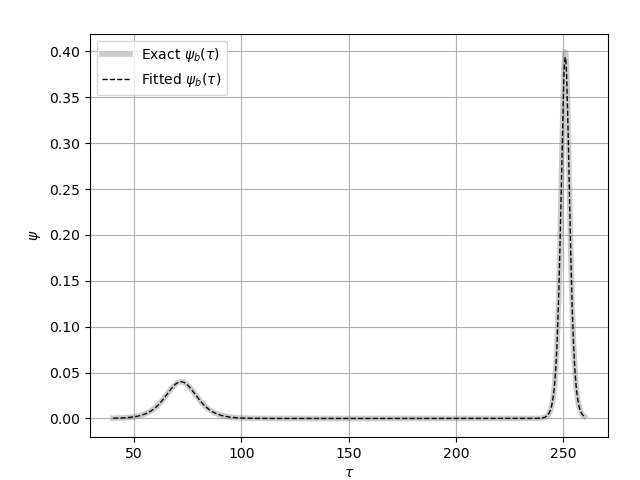}
        \caption{Fit of the recovered $\psi_{\text{b}}$ compared with the original $\psi_{\text{b}}$. Cropping was at $t = 50$ and $t = 280$.}
        \label{fig:Fitting_010_036_far}
    \end{subfigure}
    \caption{Numerical experiment with soliton parameters $q_1 = 0.1$, $q_2 = \sqrt{0.1}$, $t_1 = 75$, and $t_2 = 250$.}\label{fig:stitch_exp_2}
\end{figure}

\begin{figure}
    \begin{subfigure}{0.45\textwidth}
        \centering
            \includegraphics[width=\linewidth]{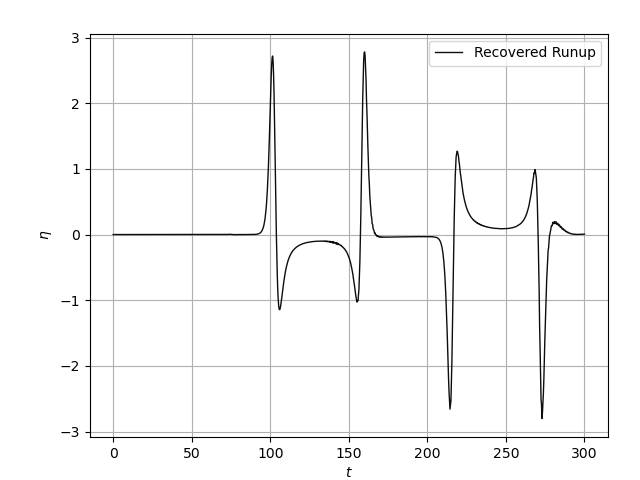}
        \caption{Run-up generated by a two-soliton solution.}
        \label{fig:Runup_036_010_far}
    \end{subfigure}
    \hfill
    \begin{subfigure}{0.45\textwidth}
        \centering
            \includegraphics[width=\linewidth]{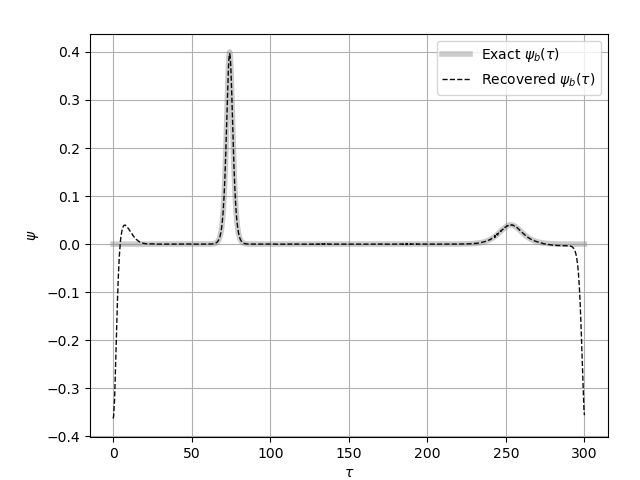}
        \caption{Recovered $\psi_{\text{b}}$ from the run-up.}
        \label{fig:Recovery_036_010_far}
    \end{subfigure}
    \\
    \centering
    \begin{subfigure}{0.45\textwidth}
            \includegraphics[width=\linewidth]{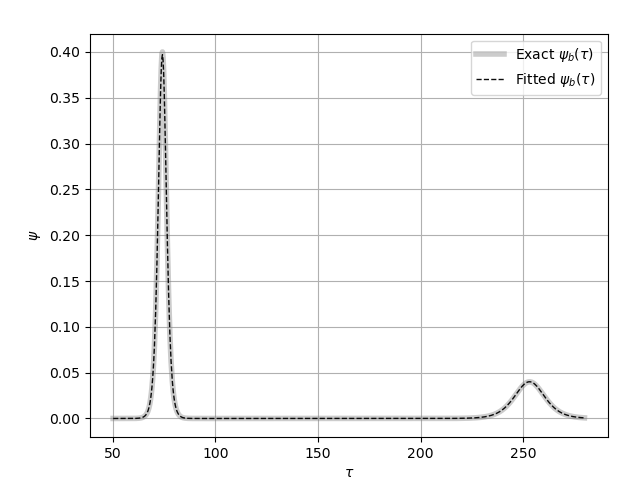}
        \caption{Fit of the recovered $\psi_{\text{b}}$ compared with the original $\psi_{\text{b}}$. Cropping was at $t = 50$ and $t = 280$.}
        \label{fig:Fitting_036_010_far}
    \end{subfigure}
    \caption{Numerical experiment with soliton parameters $q_1 = 0.1$, $q_2 = \sqrt{0.1}$, $t_1 = 75$, and $t_2 = 250$.}\label{fig:stitch_exp_3}
\end{figure}

We begin with the 2-soliton solution for \eqref{eq:boussinesq} found in \cite{hirota73} which we use as the boundary condition for $\eta(L, t)$, that is, we set
\begin{equation} \label{2soliton}
     \eta(L, t) = \frac{q_1^2 \operatorname{sech}^2 \xi_1 + q_2 ^2 \operatorname{sech}^2 \xi_2 + A \operatorname{sech}^2 \xi_1 \operatorname{sech} ^2 \xi_2}{\left( \cosh \frac{\phi}{2} + \sinh \frac{\phi}{2}\tanh \xi_1 \tanh \xi_2 \right)^2},
\end{equation}
where
\begin{equation} \label{xi}
    \xi_i = q_i L - w_i (t-t_i), \quad w_i = \epsilon_i q_i v_i,\quad v_i = \sqrt{1+4 q_i^2},
\end{equation}
\begin{equation}
    \exp(2 \phi) = \frac{(\epsilon_1 v_1 - \epsilon_2 v_2 )^2 + 12 (q_1 -q_2 )^2}{(\epsilon_1 v_1 - \epsilon_2 v_2 )^2 + 12 (q_1 +q_2 )^2},
\end{equation}
and
\begin{equation}
    A = \sinh \frac{\phi}{2} \left( (q_1 ^2 + q_2 ^2)\sinh \frac{\phi}{2} + 2 q_1 q_2 \cosh \frac{\phi}{2} \right).
\end{equation}
We then run three numerical experiments with different values of parameters \(t_i, q_i\) (Figures \ref{fig:stitch_exp_1}, \ref{fig:stitch_exp_2}, and \ref{fig:stitch_exp_3}).
Using sufficiently large $L$, gives $L \approx \sigma_L$. Therefore, $\eta(L,t)\approx \psi(\sigma_L, \tau +u)$. (Note that $L=200$, for $\alpha = 1$ and characteristic height of 5m, is 1 km).
From there, we employ the forward algorithm \eqref{eq:shoreline} in the hodograph plane to recover the corresponding shoreline conditions and, from \eqref{eq:CG_at_shore}, run-up and run-down behaviour in Figures \ref{fig:Runup_010_036_close}, \ref{fig:Runup_010_036_far}, and \ref{fig:Runup_036_010_far}.

With these shoreline conditions, we use our inversion algorithm \eqref{eq:inv_lap} to recover $\psi_{\text{b}}(\tau)$ and thus $\eta(L,t)$. Due to the bouncing waves, there are undesirable edge behaviours: Figures \ref{fig:Recovery_010_036_close}, \ref{fig:Recovery_010_036_far}, and \ref{fig:Recovery_036_010_far}. To address this, we cut off the edges of the interval. This filtering results in just the initial solitons. We then use a least squares technique to recover the soliton characteristics at the boundary: Figures \ref{fig:Fitting_010_036_close}, \ref{fig:Fitting_010_036_far}, and \ref{fig:Fitting_036_010_far}.

Note $L = 200$ for all of the following simulations.

Notice that in Figure \ref{fig:Recovery_010_036_close} the edge effect on the left is a manifestation of the original soliton. This is a demonstration of how the edge effects are related to the bouncing waves.

From \eqref{2soliton} (and general knowledge of solitons), it is known that the peak location in time and space is given when $\xi_i = 0$. Using this, we can derive the speed of each soliton. By substituting \eqref{xi} in for $\xi$ and solving for $x$, we find $x_i = -\sqrt{1+4q_i^2} \hspace{3pt}(t-t_i)$ where $x_i$ is the location of each soliton. This relationship, with the known time of the event, gives the positions of the solitons when the tsunamigenic event occurred.

\begin{figure}
    \centering
    \includegraphics[width=0.5\linewidth]{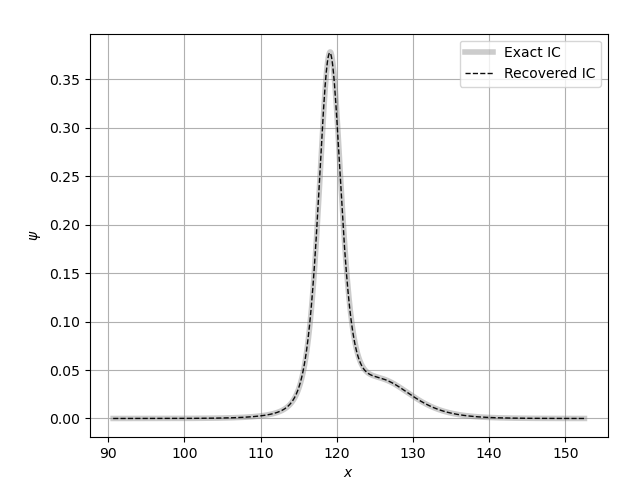}
    \caption{Recovered and exact initial conditions. Exact peaks are $x_1 = 119.55$ and $x_2 = 123.63$ and recovered peaks are $x_1 = 119.50$ and $x_2 = 123.64$.}
    \label{fig:End_036_010_far}
\end{figure}

We perform this process for the last simulation, where $q_1 = 0.1$, $q_2 = \sqrt{0.1}$, $t_1 = 75$, and $t_2 = 250$. Let's say the earthquake or tsunamigenic event happened at $t_{\text{event}} = -100$s. What we recover can be seen in Figure~\ref{fig:End_036_010_far}.

\section{Linear Shallow Water Equations}\label{sec:compare}
As another example using the piecewise nature of our model, we invert the problem proposed in \cite{Synolakis87}. This involves the same bathymetry, but SWE governs both regions. Following \cite{Synolakis87}, we use the linear SWE, which reduces to the wave equation for the free surface displacement $\eta$. Again, we assume that $\psi_{\text{b}} \approx \eta_{\text{b}}$. We also assume that the wave seen at the buoy is a travelling wave and therefore has not changed shape since it was generated. We can then do all the same steps as before in Section \ref{subsec:bouss_implementation}. Below is an implementation of this process using a soliton and an N-wave as the wave at the boundary.
\subsection{Inversion Implementation}
\begin{figure}
    \begin{subfigure}{0.45\textwidth}
        \centering
            \includegraphics[width=\linewidth]{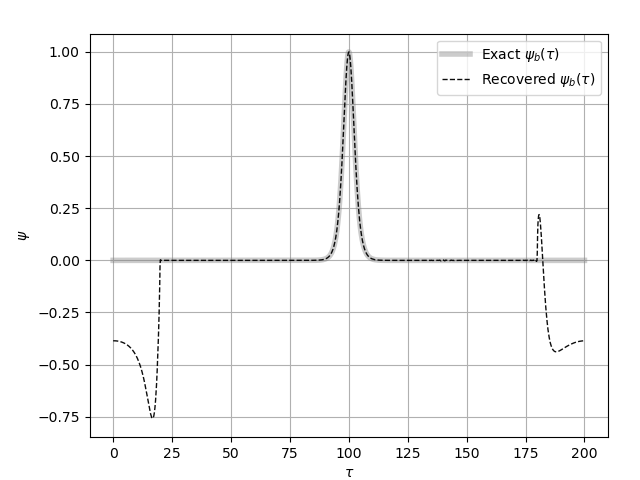}
        \caption{Recovered $\psi_{\text{b}}$ compared with input soliton.}
        \label{fig:SynolakisSoliton}
    \end{subfigure}
    \hfill
    \begin{subfigure}{0.45\textwidth}
        \centering
            \includegraphics[width=\linewidth]{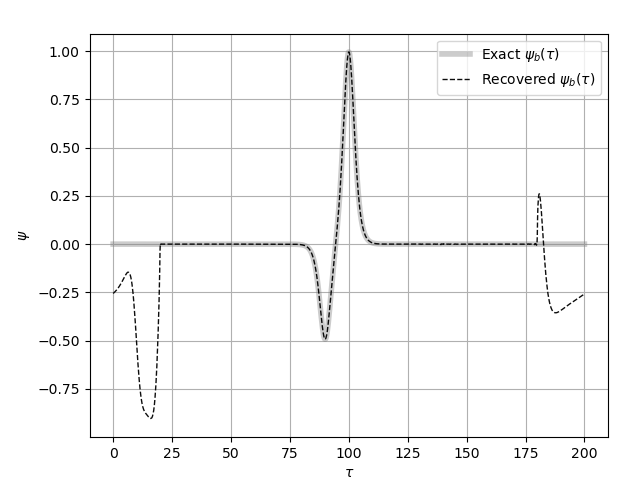}
        \caption{Recovered $\psi_{\text{b}}$ compared with input N-wave.}
        \label{fig:SynolakisNwave}
    \end{subfigure}
    \caption{Recovered and exact boundary data.}\label{fig:Synol_N_plus_sol}
\end{figure}

\begin{figure}
    \begin{subfigure}{0.45\textwidth}
        \centering
        \includegraphics[width=\linewidth]{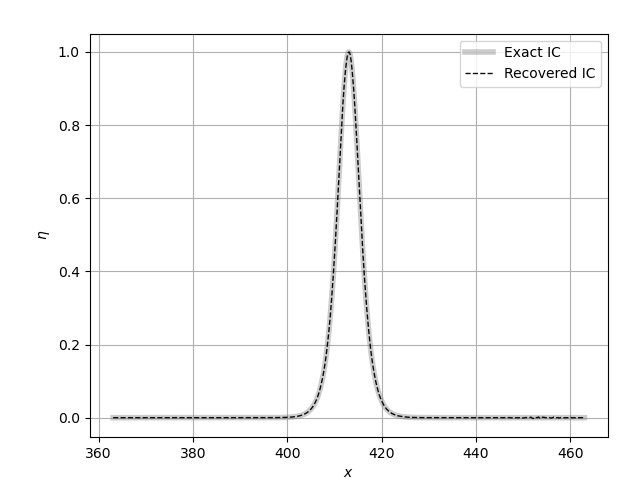}
        \caption{Recovered and exact initial conditions. Exact peak is ${x_{\text{peak}} = 412.98}$, and recovered one is ${x_{\text{peak}} = 413.12}$.\\~}
        \label{fig:FinalSynolakis}
    \end{subfigure}
    \hfill
    \begin{subfigure}{0.45\textwidth}
        \centering
        \includegraphics[width=\linewidth]{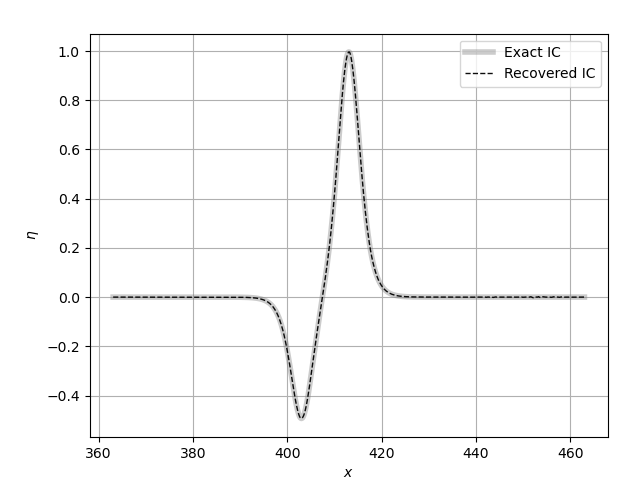}
        \caption{Recovered and exact initial conditions. Exact peak and trough are ${x_{\text{peak}} = 413.12}$ and ${x_{\text{trough}} = 402.98}$. Recovered ones are ${x_{\text{peak}} = 413.12}$ and ${x_{\text{trough}} = 402.98}$.}
        \label{fig:FinalNwaveSynolakis}
    \end{subfigure}
    \caption{Recovered and exact initial conditions.}\label{fig:final_synol}
\end{figure}

We consider two different boundary conditions
\begin{equation*}
\begin{aligned}
     \eta_{\text{b}} &= \operatorname{sech}^2((L+c(t-100))/g), \\ 
     \eta_{\text{b}} &= \operatorname{sech}^2((L+c(t-100))/g) - 0.5\operatorname{sech}^2((L+c(t-90))/g),
\end{aligned}
\end{equation*}
where $c = \sqrt{gH_0}$ and take $L = 200$. Recall that we can assume $t \approx \tau$ because we assume small $u$. Then, running the Region 1 problem (NSWE) in the forward and inverse directions returns boundary data give in Figure \ref{fig:Synol_N_plus_sol}.
After cropping, we can calculate the peak locations at the time of the event using the wave speed \(c\). We say the event happened 100 seconds before the buoy's time zero. We then plot the waves (Figure \ref{fig:final_synol}).
\subsection{Comparison}
\begin{figure}
    \centering
    \includegraphics[width=0.75\linewidth]{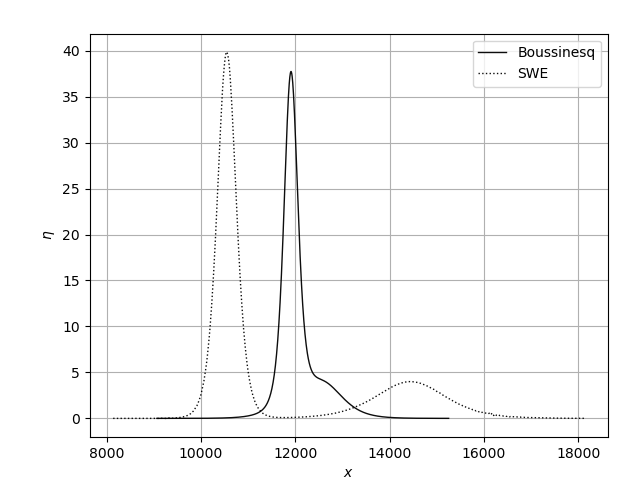}
    \caption{Comparison of initial states.}
    \label{fig:SynoVsBoussi}
\end{figure}

We then compare the initial states found by the two models using $t_{\text{event}} = -100$s and take the same \(\eta(L,t)\) as in \eqref{2soliton}.
The comparison of the Boussinesq model and SWE can be seen in Figure \ref{fig:SynoVsBoussi}.
\section{Conclusions}\label{sec:conclusions}
We have developed a robust analytical method for recovering the wave field at the virtual buoy using shoreline data. This method can be used to improve tsunami mitigation by understanding the optimal location of placing mareographs, such as DART system buoys. Numerical experiments verify the validity of our model. It is important to note that since we consider a one-dimensional model, the computational time is small, which allows for fast calculations of characteristics of fast-travelling tsunami waves. We further suggested a way to treat more complex piece-wise sloping bays, as well as to incorporate dispersion in tsunami wave propagation. One potential generalisation is to consider multi-soliton solutions for modelling tsunami propagation over the flat region, or improve accuracy by considering the boundary value problem for Boussinesq or KdV-KdV (similar to Korteweg-De Vries, but allows for multi-directional propagation of waves) equations on the half-line using techniques of \cite{Fokas}.

Improving the physicality of the bathymetry remains for future iterations. One potential course is solving the inverse problem for the multi-sloping bathymetry, similar to the one presented in \cite{Kanoglu97} within the NSWE framework.

Another important goal is to describe the set of admissible shoreline data $R(t)$, that is, to describe which run-up functions correspond to non-breaking waves.
It was shown by \cite{Rybkin14}, that a wave breaks (gradient catastrophe happens) if and only if the Jacobian of the Carrier-Greenspan transformation (or its inverse) vanishes: 
\begin{equation}
     \det \frac{\partial (\sigma,\tau)}{\partial(x,t)} = 0 \text{ or }
     \det \frac{\partial(x,t)}{\partial(\sigma,\tau)} = 0.
\end{equation}
Given run-up data $R(t)$, one can formally recover the boundary conditions, run the forward simulation, and sample these Jacobians over the computational domain to see if the wave breaks. However, we hope to find a more explicit way to describe admissible shoreline data. Another question of interest is the applicability of the iFFT; recall that we require all poles of the $\psi_{\text{sh}}$ to be to the left of the imaginary axis to implement the inverse Laplace transform using the iFFT.

Finally, we thank the anonymous referee for bringing our attention to the run-up amplification due to resonance phenomenon \cite{reso1,reso2}. We hope to investigate its impact on the inverse problem in the future.

\appendix
\section{Justification of swapping integration and summation} \label{subsec:remarks}
In \eqref{eq:shoreline} we take the Laplace transform of both sides and then switch integration and summation. Generally, integration and summation cannot be interchanged. For that, one needs the uniform (in $\tau$) convergence of the series
\begin{equation}
    \sum_{n=1}^\infty
   \frac{b_n}{\sqrt{a_n}} \int_0^ \tau \sin \left( \sqrt{a_n} (\tau- \xi) \right) \psi_{\text{b}}''(\xi) \,\d \xi.
\end{equation}
Recall that
for large \(n\), the \(n\)-th root of \(J_{0}\) is of order \(n\), so we deduce that \(a_n \sim n^2\); also \(|b_n| \sim n^{-1 / 2}\). So we deduce
\begin{equation}
     \abs{ \frac{b_n}{\sqrt{a_n}}}\sim n^{- 1/2  - 1}  = n^{- 3 /2 }.
\end{equation}
 Therefore, if we assume that \(\psi_{\text{b}}''\) is bounded (has finite \(L^{\infty}\)-norm / \(\sup\)-norm / uniform convergence norm) and restrict our attention to \(\tau \in [0, T]\) for some large \(T\) (which are reasonable assumptions from the physical point of view), we observe that the series in \eqref{eq:shoreline} converges uniformly, and so we justify interchanging integration and summation.
\bibliographystyle{apalike}
\bibliography{text.bib}
\end{document}